%% file: main.tex
\documentclass[12pt,a4paper,reqno,oneside]{amsart}
\usepackage{graphicx} 
\input{preamble}

\title{Morita equivalence and the \texorpdfstring{$\E_\infty$}{E-infinity}-Spin orientation}
\author[D.~Aretz]{David Aretz}
\address{Max-Planck-Institut f\"ur Mathematik, Vivatsgasse 7, 53111 Bonn, Germany}
\email{aretz@mpim-bonn.mpg.de}

\begin{document}
\begin{abstract}
    We exhibit a new construction of the $\E_\infty$-orientation $M\Spin \to KO$ by equating the $KO$-linear $J$-homomorphism 
    with the graded $K$-theory of Clifford algebras.
    We characterize spin structures as Morita trivializations of the Clifford algebra bundle.
    More generally, we introduce tangential structures built from bimodules.
    This allows us to deduce precise conditions on the existence of irreducible Clifford bundles, contradicting the folklore identification with spin structures.
    We also provide a homotopical interpretation of the Clifford functor as a Postnikov truncation and a higher analogue of the determinant.
    Finally, we obtain twisted variants of the spin orientation for each bimodule tangential structure.
\end{abstract}
\maketitle

\section{Introduction}
\subsection{Rapid Summary}
The main idea of this article is to unite two points of view:
\begin{enumerate}
    \item[(a)] A $\Spin$-structure on a vector bundle $V\to X$ is a Morita trivialization of the Clifford algebra bundle $\Cl(V)\to X$.
    \item[(b)] A $KO$-orientation is a nullhomotopy of the $KO$-linearized $J$-homomorphism or, equivalently, a trivialization of the $KO$-module bundle classified by the composite 
    \begin{equation*}
        X\xrightarrow{\xi_V} BO\xrightarrow{J} \Pic(\Sp) \xrightarrow{\otimes KO} \Pic(\Mod(KO))
    \end{equation*} 
\end{enumerate}
The translation is furnished by an equivalence
\begin{equation}
    \label{eq:identification}
    J(V)\otimes KO \simeq K^{\gr}\circ \Cl(V) \,.
\end{equation}
between the $KO$-linear $J$-homomorphism and the graded $K$-theory of the stable Clifford algebra.
Thus, a $\Spin$-structure induces a $KO$-orientation.
We will leverage this basic translation and expand on it in a couple of ways.
\subsection{Linear Overview}

An orientation of a vector bundle $V\to X$ is a trivialization of the determinant line bundle $\det(V)\to X$.
Similarly, in this work, a $\Spin$-structure shall be defined as a Morita equivalence $\Cl(V)\simeq \Cl(\rk V)$ between the Clifford algebra bundle of a vector space and the trivial algebra bundle $\Cl_{\rk V}$.
We take this as an excuse to develop a general theory of \emph{bimodule tangential structures}. 
After reviewing some basic definitions, Section~\ref{sec:bimodulebundles} introduces stable tangential structures $BH(A;M)\to BO$ for pairs of a finite-dimensional super algebra $A$ and a finitely generated projective module $M$.
We show that maps $X\to BH_n(A;M)$ classify a rank $n$ vector bundle $V\to X$ together with an $A$-$\Cl(V\!-\!n)$-bimodule bundle $\Sigma$ whose fibers are of \emph{type $M$}.
\begin{introtheorem}[{Theorem~\ref{thm:typedstructures}}]
    \label{thm:A}
    For $n\geq1$, the classifying space for vector bundles of rank $n$ with $(A,M)$-structure is $BH_n(A;M)$.
    The structure group is the Lie group
    \begin{equation}
        \begin{aligned}
        H_n(A;M)&=\bigl(\Pin_n^+\grotimes \sAut_A(M)^{\op}\bigr)_{\ev}\,,
        \end{aligned}
    \end{equation}
    with its natural projection map  to $O_n$.
    $BH(A;M)\to BO$ is a stable tangential structure.
\end{introtheorem}
In Section~\ref{sec:equivalencesandirreps}, we consider examples. 
Proposition~\ref{prop:spinstructures} then concludes that $\Spin$-structures are Morita equivalences $\Sigma:\Cl(V)\simeq \Cl_{\rk V}$ by setting $A=M=\R$.
Furthermore, we can recover all the groups of Freed--Hopkins' tenfold way~\cite{FreedHopkins} by looking at super division algebras.
\begin{equation*}
    \begin{array}{c|c c c c c c c c|c c}
    s&0&-1&-2&-3&4&3&2&1&0^c&1^c\\ \hline
    D(s)&\R&\Cl_{-1}&\Cl_{-2}&\Cl_{-3}&\mathbb H&\Cl_{+3}&\Cl_{+2}&\Cl_{+1}&\C&\Cxl_{-1}\\
    H(s)&\Spin&\Pin^+&\Pin^{\tilde c+}&\Pin^{h-}&\Spin^h&\Pin^{h+}&\Pin^{\tilde c-}&\Pin^-&\Spin^c&\Pin^c
    \end{array}
\end{equation*}

A popular wrong claim is that a $\Spin$-structure on $V$ is equivalent to an irreducible $\Cl(V)$ module bundle.
Using Theorem~\ref{thm:A}, we can refute this claim in dimensions $\rk(V) \not\equiv 0 \mod 8$.
This settles the question of existence of irreducible Clifford bundles completely.
\begin{proposition}[Proposition~\ref{prop:irreduciblecliffordbundles}]
    Let $V$ be a real vector bundle of rank $n$ with inner product.
    \begin{enumerate}
        \item A graded irreducible $\Cl(V)$-module bundle exists
        if and only if $V$ admits an $H(-n)$-structure.
        \item An ungraded irreducible $\lvert\Cl(V)\rvert$-module bundle exists
        if and only if $V$ admits an $H(1-n)$-structure.
    \end{enumerate}
    The indices are read modulo $8$.
    For complex Clifford modules the same statements hold with
    $H^c$ and indices modulo $2$.
\end{proposition}

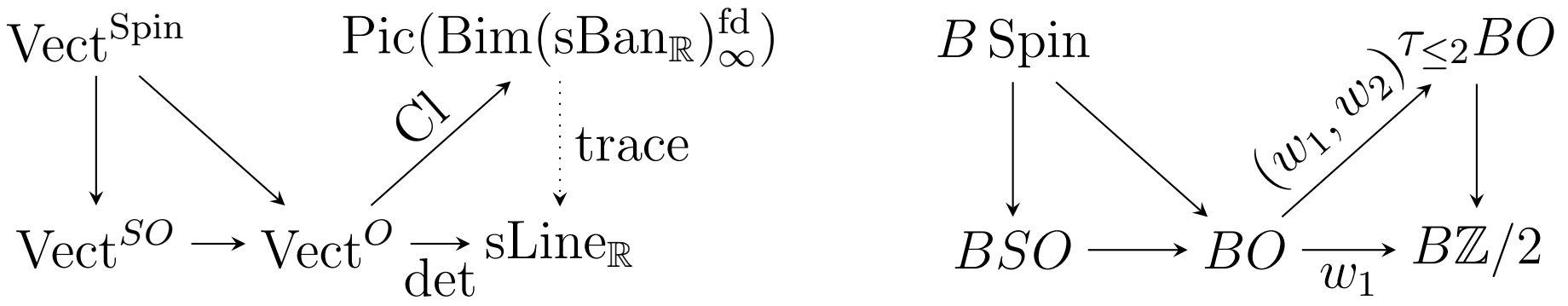
\begin{figure}[htbp]
    \centering
    \begin{tikzpicture}[>=stealth, every node/.style={inner sep=3pt}]
        \node (vspin) at (0,1.6) {$\Vect^{\Spin}$};
        \node (pic) at (3.6,1.6) {$\Pic(\Bim(\sBan_\R)^{\mathrm{fd}}_\infty)$};
        \node (vso) at (0,0) {$\Vect^{SO}$};
        \node (vo) at (1.8,0) {$\Vect^O$};
        \node (sline) at (3.6,0) {$\sLine_{\R}$};
        \draw[->] (vspin) -- (vso);
        \draw[->] (vspin) -- (vo);
        \draw[->] (vso) -- (vo);
        \draw[->] (vo) -- node[below] {$\det$} (sline);
        \draw[->] (vo) -- node[midway, sloped, above] {$\Cl$} (pic);
        \draw[->,dotted] (pic) -- node[right] {$\operatorname{trace}$} (sline);
    \end{tikzpicture}%
    \hspace{1cm}%
    \begin{tikzpicture}[>=stealth, every node/.style={inner sep=3pt}]
        \node (bspin) at (0,1.6) {$B\Spin$};
        \node (postnikov) at (3.6,1.6) {$\tau_{\leq 2}BO$};
        \node (bso) at (0,0) {$BSO$};
        \node (bo) at (1.8,0) {$BO$};
        \node (bz2) at (3.6,0) {$B\Z/2$};
        \draw[->] (bspin) -- (bso);
        \draw[->] (bspin) -- (bo);
        \draw[->] (bso) -- (bo);
        \draw[->] (bo) -- node[below] {$w_1$} (bz2);
        \draw[->] (bo) -- node[midway, sloped, above] {$(w_1,w_2)$}(postnikov);
        \draw[->] (postnikov) -- (bz2);
    \end{tikzpicture}
    \caption{The categorical Postnikov/Whitehead tower (left) and its stable homotopical
    counterpart (right), with the corresponding Postnikov stages in the right-hand column of each.}
    \label{fig:categorical-whitehead-tower}
\end{figure}

In Section~\ref{sec:categoricalpostnikov}, we provide yet another, homotopy-theoretical, proof that $\Spin$-structures are Morita trivializations by showing that the
Clifford functor can be interpreted as a categorical Postnikov truncation~(Figure~\ref{fig:categorical-whitehead-tower}).
In particular, the stable Clifford functor fits into a fiber sequence of $\E_\infty$-spaces:
\begin{equation}
    \label{eq:fibseq1}
    B\Spin \to BO \xrightarrow{\Cl} \Pic(\Bim(\sBan_\R)^{\mathrm{fd}}_\infty) \,,
\end{equation}
with target the Morita $(\infty,1)$-category of super Banach algebras, which is recalled in Section~\ref{sec:picardgroupoid}.

Section~\ref{sec:naturaltransformation} recalls briefly the construction of graded $K$-theory $K^{\gr}$ from \cite{aretzstehouwer} and a model of the $J$-homomorphism to prove the following theorem.
\begin{introtheorem}[{Theorem~\ref{thm:spin-orientation-proof}}]
    \label{thm:B}
    There is an invertible lax symmetric monoidal natural transformation making the following diagram of lax symmetric monoidal functors commute:
    \begin{equation*}
        \begin{tikzcd}
            \mathcal{I}nn \ar[r,"J"]\ar[d,"\Cl"] &\Sp \ar[d,"\otimes KO"]\\
        \Bim(\sBan_\R)_\infty \ar[r,"K^{\gr}"] & \Mod(KO)
        \end{tikzcd}
    \end{equation*}
    Here, $\Inn$ is the $\infty$-category of inner product spaces.
\end{introtheorem}
We will show that this has far-reaching consequences. 
Using Theorem~\ref{thm:B}, the fiber sequence~\eqref{eq:fibseq1} provides a nullhomotopy of the $KO$-linear $J$-homomorphism.
By \cite{AndoBlumbergGepner}, recalled in Section~\ref{sec:thomspectra}, such nullhomotopies precisely correspond to multiplicative orientations.
That is:
\begin{corollary}
    There is a map of $\E_\infty$-rings from the spin Thom spectrum to real topological $K$-theory 
    \begin{equation*}
        M\Spin \longrightarrow KO \,.
    \end{equation*}
\end{corollary}
We shall emphasize here that this construction is a) constructive and b) does not use functional analysis.
It is a homotopy-coherent version of the Atiyah--Bott--Shapiro construction~\cite{ABS} and coincides with the refinements of \cite{AndoHopkinsRezk,JoachimCoherences} as shown in Appendix~\ref{app:comparison}.

We furthermore prove that we can construct a lot more $K$-theory orientations for the bimodule tangential structures from Section~\ref{sec:bimodulebundles} in a multiplicative manner.
\begin{introtheorem}[{Theorem~\ref{thm:multiplicativeorientations}}]
    \label{thm:C}
    Let $A$ be a super Banach algebra and let $M$ be a finitely generated projective $A$-module.
     There are maps of spectra
                \begin{equation}
                    \label{eq:twistedorientation1}
                    MH(A;M) \to K^{\gr}_A \,,
                \end{equation}
    satisfying the following:
    \begin{enumerate}[label=(\roman*), ref=(\roman*)]
        \item\label{item:C-i} Let $MH(A):=\bigoplus_{[M]} MH(A;M)$, where $[M]$ runs over types of finitely generated projective $A$-modules.
                     $MH$ is the left Kan extension of the $J$-homomorphism along the stable Clifford functor $\Cl$.
        \item\label{item:C-ii} $MH$ defines a lax symmetric monoidal functor $MH:\Bim(\sBan_\R)_\infty\to \Sp$
                     and there is a universal lax symmetric monoidal natural transformation $MH\Rightarrow K^{\gr}$.
    \end{enumerate}
    Therefore, the maps \eqref{eq:twistedorientation1} are multiplicative, i.e.\ there are coherently commutative squares
    \begin{equation*}
        \begin{tikzcd}
            MH(A;M)\otimes MH(B;N) \ar[r]\ar[d]& K^{\gr}_{A}\otimes K^{\gr}_B \ar[d]\\
            MH(A\grotimes B;M\grotimes N) \ar[r]& K^{\gr}_{A\grotimes B}
        \end{tikzcd}
    \end{equation*}
    In particular, inserting $A=M=\R$ shows that \eqref{eq:twistedorientation1} are $M\Spin\to KO$-module maps.
    \begin{enumerate}[label=(\roman*), ref=(\roman*), start=3]
        \item\label{item:C-iii} Let $1\le n\le\infty$ and let $A$ be an $\E_n$-algebra in $\Bim(\sBan_\R)_\infty$.
            Then $MH(A)\to K^{\gr}_A$ refines to a map of $\E_n$-rings.
            E.g.\ if $A$ is a graded-commutative super Banach algebra and $n=\infty$,
            then $MH(A;A)$ is a subring of $MH(A)$ and $MH(A;A)\to K^{\gr}_A$ is an $\E_\infty$-$H(A;A)$-orientation.
        \item\label{item:C-iv} If $K^{\gr}_A$ is an $\E_1$-ring, e.g.\ as in~\ref{item:C-iii}, then for every $M$
            the map~\eqref{eq:twistedorientation1} is an $H(A;M)$-orientation of $K^{\gr}_A$
            if and only if $[M]\in\pi_0K^{\gr}_A$ is a unit.
    \end{enumerate}
\end{introtheorem}
Section~\ref{sec:examples} finally gathers examples: we can recover the tenfold way of Freed--Hopkins~\cite{FreedHopkins} and point out partial multiplicativity.
We obtain multiplicative orientations $MH^G\to K^{\gr}_{\R^f[G]}$ for 
finite fermionic symmetry groups $G$ and their spacetime structure group $H^G=(\Pin_n^+\grotimes G)_{ev}$.
For infinite discrete groups we obtain a real Bost assembly map $M\Spin \wedge B\Gamma_+\to K_{\ell^1(\Gamma)}$.
Finally, we also sketch the relation to the Sullivan orientation $MSO\to KO[\tfrac{1}{2}]$ as a map of $\E_\infty$-rings.

\subsection{Related Work}
\label{sec:relatedwork}
Our description of spin structures as Morita equivalences was inspired by Stolz and Teichner~\cite{StolzTeichnerElliptic}, who define a spin structure on $V$ as a bundle of irreducible $\Cl(V)$-$\Cl_{\rk V}$-bimodules.
The spacetime structure groups of Freed and Hopkins~\cite{FreedHopkins} motivated our study of super division algebras.
The $\E_\infty$-refinement $M\Spin\to KO$ of the Atiyah--Bott--Shapiro orientation~\cite{ABS} was previously obtained by obstruction theory~\cite{AndoHopkinsRezk} and from point-set models of $KO$ built from Fredholm operators~\cite{JoachimCoherences,HebestreitJoachim}; our construction uses neither.
Ebert~\cite{EbertIndex} constructs maps $MT\theta_A\to K_A$ for unstable tangential structures for super $C^*$-algebras using Fredholm models; Theorem~\ref{thm:C} provides stable, multiplicative analogues.
Kamel and Mooers~\cite{kamelmooers} also study a generalized spin orientation of the $K$-theory of $C^*$-algebras in work in progress.
\subsection{Acknowledgements}
The author thanks Christian Blohmann, Luuk Stehouwer and Peter Teichner for helpful input, feedback and encouragement.
Further the author thanks Johannes Ebert for interesting discussions and Yigal Kamel and Fredrick Mooers for sharing openly their work-in-progress.

The author acknowledges the hospitality of the Max Planck Institute for Mathematics in Bonn.

\subsection{AI usage disclosure}
An earlier version of this article was written without AI assistance.
In revising it, the author used the language models Opus~5.5 and Astra~6 to check and explore proofs, find errors, draft corrections and search the literature; they also suggested Section~\ref{sec:sullivan}.

\subsection{Conventions}
We freely use the language of $\infty$-categories~\cite{HTT,HA}.
$\Spc$ denotes the $\infty$-category of spaces, $\Sp$ the $\infty$-category of spectra.
All concepts are by default homotopical, e.g.\ $\colim$ refers to the $(\infty,1)$-categorical colimit.

\section{Bimodule Bundles}
\label{sec:bimodulebundles}

\subsection{Super Algebras}
We recall the Morita theory for finite-dimensional super algebras.
This material is mostly well-known or folklore, we refer to \cite{aretzstehouwer} for more details.
A super vector space is a $\Z/2$-graded $\mathbb{F}$-vector space $V=V_0\oplus V_1$.
We will work mainly with $\mathbb{F}=\R$.
Homomorphisms of super vector spaces are required to preserve the grading.
The category of super vector spaces $\sVect_{\mathbb{F}}$ admits a symmetric monoidal structure $\grotimes$, whose braiding is the Koszul braiding
\[
    \beta:V\grotimes W\to W\grotimes V\,, \qquad v\otimes w\mapsto (-1)^{|v||w|}w\otimes v\,.
\]

A super algebra is an algebra object in $(\sVect_{\mathbb{F}},\grotimes)$.
A (graded left) module over a super algebra $A$ is an $A$-module object $M$ in $(\sVect_{\mathbb{F}},\grotimes)$.
Explicitly, this requires an action $A\grotimes M\to M$ that is unital and associative in the appropriate sense.
We write $\Mod^{C_2}_A$ for the category of finitely generated projective left modules and grading-preserving homomorphisms.
An $A$-module $M$ is \emph{finitely generated projective} if it is a direct summand of $A\grotimes \R^{p|q}$.
$M$ is \emph{irreducible} if it does not contain a proper graded $A$-submodule.
A super algebra is a \emph{super division algebra} if every nonzero homogeneous element $a\in A_0\cup A_1$ is invertible.

\begin{lemma}[Schur's Lemma]
    Let $M$ be an irreducible $A$-module. 
    Then the super algebra $\mathrm{sEnd}_A(M)$ of endomorphisms of $M$ which graded-commute with $A$ is a super division algebra over $\mathbb{F}$.
\end{lemma}
\begin{remark}
    \label{rmk:handedness}
    The super-opposite
    $A^{\op}$ has multiplication
    $a\cdot_{\op}b=(-1)^{|a||b|}ba$.
    We have $\mathrm{sEnd}_A(A)\cong A^{\op}$.
\end{remark}
\begin{lemma}
    There is an equivalence between
    \begin{enumerate}
        \item $\sVect$-linear, $\oplus$-preserving functors $\Mod^{C_2}_A\to \Mod_{B}^{C_2}$.
        \item $B$-$A$-bimodules $M$ that are finitely generated projective as $B$-modules.
    \end{enumerate}
    This is induced by mapping $M\mapsto M\grotimes_A -$
\end{lemma}
We write $M:A\squigrightarrow B$, when viewing a $B$-$A$-bimodule as a morphism from $A$ to $B$.
There is a symmetric monoidal Morita bicategory $(\Mor(\sVect),\grotimes)$ of super algebras, bimodules and grading-preserving bimodule homomorphisms.
A \emph{Morita equivalence} is an isomorphism in $\Mor(\sVect)$ or, equivalently, an $\sVect$-linear additive equivalence of module categories.
\begin{definition}
    Let $V$ be a finite-dimensional real vector space with a nondegenerate symmetric bilinear form $B$.
    The \emph{Clifford algebra} $\Cl(V,B)$ is defined as the free (noncommutative) algebra on the vector space $V$ modulo the relations $vw + wv = 2 B(v,w)$ for all $v,w\in V$.
    $\Cl(V)$ has a natural $\Z/2$-grading where all generators $v\in V$ are odd.
    
    We write $\Cl_n=\Cl_{+n}$ for $\Cl(\R^n,\langle,\rangle)$ and $\Cl_{-n}$ for $\Cl(\R^n,-\langle,\rangle)$, where $\langle,\rangle$ is the standard Euclidean inner product.
    We write $\Cxl_n = \Cl_n\otimes \C$.
\end{definition}
Under the super tensor product there is a natural isomorphism of super algebras
\begin{equation*}
    \Cl(V,B)\grotimes \Cl(W,B')\cong \Cl(V\oplus W, B\oplus B')\,, \qquad \Cl(V,B)^{\op}\cong \Cl(V,-B) \,.
\end{equation*}
If $A$ is a super algebra, we write $|A|$ for the underlying ungraded algebra.
\begin{example}
    We have $\lvert\Cl_0\rvert=\R$, $\lvert\Cl_{-1}\rvert=\C$, $\lvert\Cl_{-2}\rvert=\mathbb{H}$, $\lvert\Cl_{+1}\rvert\cong \R\oplus \R$.
\end{example}

\begin{example}
    There is a natural equivalence of categories $\Mod_{A}^{C_2}\simeq \Mod_{\lvert A\grotimes \Cl_{+1}\rvert}$ \cite[Prop.~3.36]{aretzstehouwer}.
    For the super algebra $\R$ with trivial grading this amounts to the equivalence $\sVect\cong \Mod_{\R\oplus \R}\cong \Vect^{\times 2}$.
\end{example}
\begin{example}
    $\R$ is Morita equivalent to $\Cl_{1,1}=\Cl_{+1}\grotimes \Cl_{-1}$ and furthermore $\Cl_{n,n}\simeq \R$.
    A bimodule inducing this equivalence is $\Cl_n$, regarded as a left $\Cl_n\grotimes \Cl_n^{\op}$-module.
    Another nonisomorphic choice is the parity shift $\Pi \Cl_n$. 
\end{example}
\begin{example}
    $\R$ is Morita equivalent to $\Cl_{8k}$ for all $k\in \Z$.
    This can be called \emph{algebraic Bott periodicity}.
\end{example}

Let $(V,\langle,\rangle)$ be an inner product space.
We write $\Cl(V\!-\!n)$ for $\Cl(V)\grotimes \Cl_{-n}$.
By uniqueness of inner products, the super algebra $\Cl(V\!-\!\dim V)$ is isomorphic to $\Cl_{\dim V,\dim V}$ and hence Morita equivalent to the ground field $\R$.
As $\R$ has two irreducible graded modules, it follows that $\Cl(V\!-\!\dim V)$ has two irreducible representations up to isomorphism: $E$ and its parity shift $\Pi E$.

\begin{definition}
    \label{def:type}
    Let $A$ be a super algebra and $M$ a super representation.
    An $A$-$\Cl(V\!-\!\dim V)$-bimodule $\Sigma$ is \emph{of type $M$} if
    the functor 
    \begin{equation*}
        \Sigma \grotimes_{\Cl(V\!-\!\dim V)} -: \Mod^{C_2}_{\Cl(V\!-\!\dim V)} \to \Mod^{C_2}_A \,.
    \end{equation*}
    maps the irreducible representations $\{E,\Pi E\}$ to $A$-modules which are isomorphic to $M$ or $\Pi M$.
    We say that the type $M$ is \emph{oriented} if $M$ and $\Pi M$ are not isomorphic and \emph{unoriented} if $M\cong \Pi M$.
\end{definition}
Being of the same type is an equivalence relation on $A$-modules which identifies $M$ and $\Pi M$.
\begin{example}
    \label{ex:irreducibleClnnmodules}
    Let $g\in O(\R^n)$.
    The algebra automorphism $\Cl(g):\Cl_n\to \Cl_n$ induces the bimodule $(\Cl_n)_g:\Cl_n\squigrightarrow\Cl_n$, where the right action is restricted along $\Cl(g)$.
    The isomorphism class of $(\Cl_n)_g$ depends on $\det(g)$.

    We will reason about $\Cl_{n,n}$-modules instead of $\Cl_n$-bimodules.
    The two isomorphism classes of irreducibles are distinguished by the action of the even volume element $\omega=e_1\dots e_nf_1\dots f_n$ with $\omega^2=1$ and which satisfies $\omega a=(-1)^{|a|}a\omega$.
    Hence, on every representation the composition of $\omega$ with the grading operator $(-1)^F$ defines an even operator which commutes with the $\Cl_{n,n}$-action. 
    If the representation is irreducible, by Schur's Lemma $\omega$ must act as $\pm (-1)^F$.
    The sign distinguishes the two irreducible representations.
    On $(\Cl_n)_g$ the volume element acts as $e_1\dots e_n\,g(f_1)\dots g(f_n)=\det(g)\omega$.
    Hence, $(\Cl_n)_g\cong \Pi^{\epsilon(g)} \Cl_n$ where $\epsilon(g)=(1-\det(g))/2\in\Z/2$.
\end{example}
\begin{lemma}
    \label{lem:type}
    Let $\Sigma,\Sigma'$ be $A$-$\Cl(V\!-\!\dim V)$-bimodules.
    They have the same type if and only if there are bimodule isomorphisms $\Sigma\cong \Sigma'$ or $\Sigma\cong \Pi \Sigma'$.
\end{lemma}
\begin{proof}
    We can assume $V=\R^n$.
    Suppose they have the same type.
    Let a $\Cl_{n,n}$-$\R$-bimodule $P$ and an $\R$-$\Cl_{n,n}$-bimodule $Q$ constitute a Morita equivalence $\R\simeq \Cl_{n,n}$.
    $P$ is necessarily irreducible and by assumption $\Sigma\grotimes_{\Cl_{n,n}} P\cong M$ (or $\cong \Pi M$) which is furthermore isomorphic to $\Sigma'\grotimes_{\Cl_{n,n}} P$.
    Applying $-\grotimes Q$ and the isomorphism $P\grotimes Q\cong \Cl_{n,n}$ as $\Cl_{n,n}$-bimodules, we obtain an isomorphism $\Sigma\cong \Sigma'$ (or $\cong \Pi \Sigma'$).
    The converse is trivial.
\end{proof}

\subsection{Fermionic Groups}
\begin{definition}
    A \emph{fermionic group} is a group $G$ together with a group homomorphism $|\cdot|:G\to \Z/2$, called \emph{degree}, and a central element $c\neq 1$ with $c^2=1$ and $|c|=0$.
    We can view this as a central extension:
    \begin{equation*}
        \begin{tikzcd}
            1\ar[r]& \Z/2\ar[r]&  G \ar[r] \ar[rd,"|\cdot |"]&   G_b \ar[r]\ar[d]& 1\\
            & & & \Z/2 &
        \end{tikzcd}
    \end{equation*}
    $G_b$ is called the \emph{underlying (bosonic) group} of $G$.
    $G_{\ev}:=\ker(|\cdot|)$ is the subgroup of even elements.
\end{definition}
\begin{remark}
    \label{rmk:fermionicclassification}
    It is often beneficial to view a fermionic group as a map $BG_b\to B\Z/2 \times B^2 \Z/2$, where the map to the first factor encodes $|\cdot|$ and the map to the second factor encodes the central extension.
    This explains the classification of fermionic group structures by group cohomology classes $w_1(G)\in H^1(G_b,\Z/2),w_2(G)\in H^2(G_b,\Z/2)$~\cite[Def.~3.2.16 and Rem.~3.2.18]{luukphd}.
\end{remark}
We will use that there is an equivalence of spaces $(w_1,w_2):\tau_{\leq 2}BO\to B\Z/2\times B^2\Z/2$. Homotopy-theoretically, a fermionic group is just a map
\begin{equation*}
    BG_b\longrightarrow \tau_{\leq 2}BO \,.
\end{equation*}
\begin{example} 
    Let $A$ be a finite-dimensional super algebra. 
    The set of invertible homogeneous elements in $A$, denoted $A^\times=A_0^\times\cup A_1^\times$, forms a fermionic Lie group.
    The degree is given by the grading and the central element is $c=-1$.
\end{example}
\begin{example}
    The group $\Pin_n^{\pm}\subset \Cl_{\pm n}^{\times}$ is generated by the unit vectors in $\R^n$.
    Its even subgroup is $\Spin_n$, and the twisted adjoint action gives the covering maps
    \begin{equation*}
        \Pin_n^{\pm}\longrightarrow O_n
        \quad\text{and}\quad
        \Spin_n\longrightarrow SO_n
    \end{equation*}
    with kernel $\{\pm 1\}$.
\end{example}
\begin{example}
    \label{ex:pinfermionic}
    Consider the canonical truncation map $BO_n\to \tau_{\leq 2}BO$.
    This classifies the extension $\Z/2\to \Pin_n^+\to O_n$ with nontrivial grading.
\end{example}
\begin{proposition}
    \label{prop:fiberseqeven}
    Let $G$ be a fermionic group with classifying map $BG_b\to \tau_{\leq 2}BO$ and let $G_{\ev}$ be its even subgroup.
    Suppose that $G$ contains an odd element.
    There is a fiber sequence
    \begin{equation*}
        BG_{\ev} \to BG_b \to \tau_{\leq 2}BO \,.
    \end{equation*}
\end{proposition}
\begin{proof}
    There are horizontal and vertical fiber sequences:
    \begin{equation*}
        \begin{tikzcd}
            BG_{\ev}\ar[r]\ar[d]&B(G_b)_{\ev} \ar[r]\ar[d] & B^2\Z/2 \ar[d]\\
            F\ar[r]\ar[d]&B G_b \ar[r] \ar[d,"|-|"] & \tau_{\leq 2} BO \ar[d] \\
            * \ar[r]&B\Z/2 \ar[r,equals]& B\Z/2  
        \end{tikzcd}
    \end{equation*}
    Here $(G_b)_{\ev}=G_{\ev}/\langle c\rangle$ is the kernel of the degree $G_b\to\Z/2$.
    Hence, $BG_{\ev}\simeq F$.
\end{proof}
\begin{definition}
    Let $G,H$ be fermionic groups.
    The \emph{fermionic tensor product} $G\grotimes H$ is generated by $g\grotimes h$, $g\in G$, $h\in H$ with multiplication
    \begin{align*}
        (g_1\grotimes h_1)(g_2\grotimes h_2) = c^{|g_2||h_1|} g_1g_2\grotimes h_1h_2 \,,
    \end{align*}
    subject to the relation
    \begin{equation*}
        c=c_G\grotimes 1 = 1\grotimes c_H \,.
    \end{equation*}
    The degree is $|g\grotimes h|=|g|+|h|$ and the central element is $c$.
    As a set, we have $G\grotimes H=G\times H/\langle (c_G,c_H)\rangle $.
\end{definition}
\begin{remark}
    The tensor product of fermionic groups is induced from the $\E_\infty$-structure $\oplus$ on $B\Z/2\times B^2\Z/2=\tau_{\leq 2}BO$.
    In terms of the cohomology classes of Remark~\ref{rmk:fermionicclassification}, the fermionic tensor product corresponds to $w_1(G\grotimes H)=w_1(G)+w_1(H)$ and $w_2(G\grotimes H)=w_2(G)+w_2(H)+w_1(G)\cup w_1(H)$.
\end{remark}
\begin{example}
    Consider a tensor product $A\grotimes B$ of super algebras.
    Then, there is a natural map $A^\times \grotimes B^\times \to (A\grotimes B)^\times$.
\end{example}
The \emph{opposite fermionic group} $G^{\op}$ has multiplication $g\cdot h:=c^{|g||h|}hg$.
In terms of classifying maps, this corresponds to the postcomposition with inversion $BG_b \to \tau_{\leq 2}BO \xrightarrow{-} \tau_{\leq 2}BO$. 
\begin{example}
    \label{ex:oppositefermionicgroup}
    The fermionic group $(A^{\op})^\times$ is $(A^\times)^{\op}$.
\end{example}
\begin{example}
    We have $(\Pin_n^+)^{\op}\cong \Pin_n^-$.
\end{example}

\subsection{Bimodule Bundles of Type $(A,M)$}
This section introduces a stable tangential structure built from bimodule bundles and identifies its classifying space.

A \emph{locally trivial algebra bundle} is an algebra object in vector bundles that is locally isomorphic, compatibly with multiplication, to a product $U\times B_0\to U$, where $B_0$ is an algebra.
A \emph{locally trivial module bundle} $N\to X$ over a locally trivial algebra bundle $B\to X$ is a vector bundle with a fiberwise $B$-action that admits simultaneous trivializations $B|_U\cong U\times B_0$ and $N|_U\cong U\times N_0$, where $N_0$ is a $B_0$-module and the isomorphisms are compatible with the action.

Let $V\to X$ be a vector bundle equipped with an inner product $\langle,\rangle$.
Applying the functor $\Cl$ to each fiber $(V_x,\langle,\rangle)$ we obtain the \emph{Clifford algebra bundle} $\Cl(V)\to X$.
This defines a locally trivial bundle of super algebras over $X$.
We write $\Cl(V\!-\!n)$ for the super algebra bundle $\Cl(V)\grotimes \Cl_{-n}$.

The fibers of the super algebra bundle $\Cl(V\!-\!\rk V)$ are noncanonically isomorphic to $\Cl_{\rk V,\rk V}$ and are hence Morita equivalent to $\R$.
Each fiber has two inequivalent irreducible representations which are related by the parity shift operator $\Pi$.

\begin{definition}
    Let $A$ be a finite-dimensional super algebra and let $M$ be an $A$-module.\footnote{By convention, $M$ is a super vector space and a graded representation.}
    An \emph{$(A,M)$-structure} on a vector bundle $(V,\langle, \rangle)$ with inner product is a locally trivial $A$-$\Cl(V\!-\!\rk V)$-bimodule bundle $\Sigma$ such that each fiber $\Sigma_x$ is of type $M$ (Definition~\ref{def:type}).
    Emphasizing the algebra, we also say that $\Sigma$ is a bimodule bundle \emph{of type $(A,M)$}.

    An isomorphism $(f,\phi):(V,\Sigma)\to (V',\Sigma')$ of $(A,M)$-structures is an isometric isomorphism $f: V\to V'$ together with an isomorphism of bimodule bundles
    \begin{equation*}
        \phi: \Sigma \xrightarrow{\cong} \Sigma'_f \,,
    \end{equation*}
    where $\Sigma'_f$ denotes $\Sigma'$ with the right Clifford action restricted along the algebra
    isomorphism induced by $f$.
\end{definition}
We shall think of the bimodule bundle as a morphism of super algebra bundles $\Sigma:\Cl(V\!-\!\rk V)\squigrightarrow \underline{A}$.
An isomorphism $(f,\phi)$ is then a $2$-cell $\phi:\Sigma\Rightarrow\Sigma'\circ\Cl(f)$, since precomposing with the bimodule induced by $\Cl(f)$ restricts the right action, i.e.\ $\Sigma'\circ\Cl(f)=\Sigma'_f$ (cf.\ Example~\ref{ex:irreducibleClnnmodules}):
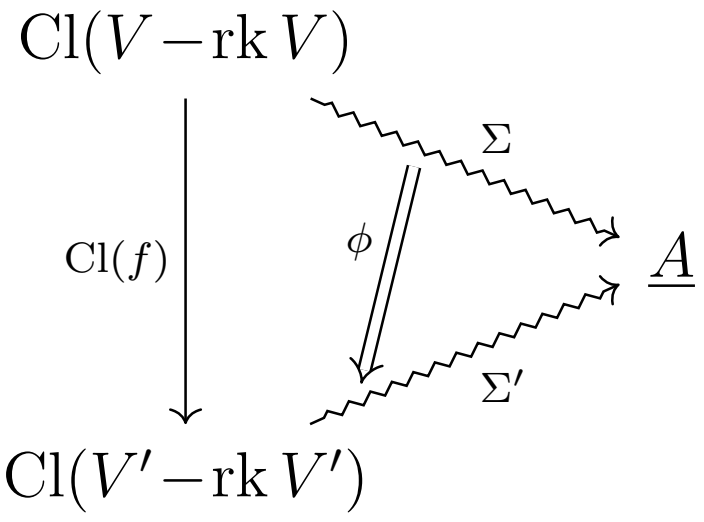
\begin{equation*}
    \begin{tikzcd}[column sep=large, row sep=normal]
    {\Cl(V\!-\!\rk V)} \\
    & \underline{A} \\
    {\Cl(V'\!-\!\rk V')}
    \arrow["\Sigma", ""{name=0, anchor=center, inner sep=0, pos=0.35}, squiggly, from=1-1, to=2-2]
    \arrow["{\Cl(f)}"', from=1-1, to=3-1]
    \arrow["{\Sigma'}"', ""{name=1, anchor=center, inner sep=0, pos=0.15}, squiggly, from=3-1, to=2-2]
    \arrow["\phi"', Rightarrow, shorten=3pt, from=0, to=1]
    \end{tikzcd}
\end{equation*}

\begin{definition}
    We write $\Vect_{(A,M)}(X)$ for the groupoid of vector bundles with $(A,M)$-structure over $X$.
    Via pullback, this defines a functor $\Vect_{(A,M)}:\Top^{\op} \to \Grpd$.
\end{definition}
\begin{example}[Local model]
    \label{ex:localmodel}
    Let $n\geq 1$.
    Regarding $\Cl_n$ as a $\Cl_n$-bimodule, i.e.\ a right $\Cl(\R^n\!-\!n)$-module, the \emph{local model} of a type $(A,M)$-structure in rank $n$ is $(\R^n,\Sigma_n)$, where
    \begin{equation}
        \label{eq:localmodel}
        \Sigma_n:=\Cl_n\grotimes M:\Cl(\R^n\!-\!n)\squigrightarrow A \,.
    \end{equation}
    This bimodule may be viewed as the composite $\Cl(\R^n\!-\!n)\simeq \R \xsquigrightarrow{M} A$ of the Morita equivalence $\Cl(\R^n\!-\!n)\simeq \R$ and the $A$-module $M$.
    Up to isomorphism, it is the only $(A,M)$-structure on an $n$-dimensional inner product space $V$.
    Indeed, after choosing an isometry $V\cong\R^n$, it suffices to show that any two $(A,M)$-structures $(V,\Sigma)$ and $(V,\Sigma')$ are isomorphic.
    By Lemma~\ref{lem:type}, there is a bimodule isomorphism $\phi:\Sigma\cong \Sigma'$ or $\phi:\Sigma\cong \Pi \Sigma'$.
    In the first case, $(\id,\phi):(V,\Sigma)\to (V,\Sigma')$ is an isomorphism.
    In the second case, pick $f\in O(V)$ with $\det(f)=-1$; Example~\ref{ex:irreducibleClnnmodules}, applied to the factor $\Cl(V)$ of $\Cl(V\!-\!n)$, gives $\Cl(V\!-\!n)_f\cong\Pi\Cl(V\!-\!n)$ and hence $\psi: \Pi \Sigma'\cong \Sigma'\circ\Cl(f)=\Sigma'_f$, and $(f,\psi\circ \phi)$ is an isomorphism.
\end{example}
Consequently, every $(A,M)$-structure on a rank $n$ vector bundle is locally isomorphic to the trivial bundle with fiber the local model $(\R^n,\Sigma_n)$.
Its transition functions take values in the closed subgroup of automorphisms of the local model
\[
    \Aut(\R^n,\Sigma_n)\subset O(n)\times \GL(\Sigma_n)\,,
\]
consisting of pairs $(f,\phi)$ such that $\phi$ is an even $A$-$\Cl_{n,n}$-bimodule isomorphism $\Sigma_n\to (\Sigma_n)_f$.

\begin{remark}
    \label{rmk:multiplicationofstructures}
    Vector bundles with bimodule tangential structure can be added: For $(V,\Sigma_V)$ in $ \Vect_{(A,M)}(X)$ and $(W,\Sigma_W)$ in $\Vect_{(B,N)}(X)$, there is
    \begin{equation}
        (V\oplus W, \Sigma_V\grotimes \Sigma_W) \in \Vect_{(A\grotimes B, M\grotimes N)}(X) \,.
    \end{equation}    
\end{remark}

\begin{theorem}
    \label{thm:typedstructures}
    Let $A$ be a finite-dimensional super algebra and $M\neq 0$ a finitely generated projective module.
    For $n\geq1$, the classifying space for vector bundles of rank $n$ with $(A,M)$-structure is $BH_n(A;M)$.
    For every paracompact locally compact Hausdorff space $X$ we have
    \[
        \Vect^{\rk=n}_{(A,M)}(X)/\text{isom.}
        \cong \pi_0[X,BH_n(A;M)].
    \]
    The structure group is
    the Lie group
    \begin{equation}
        \label{eq:grouphom}
        \begin{aligned}
        H_n(A;M)&=\bigl(\Pin_n^+\grotimes \sAut_A(M)^{\op}\bigr)_{\ev}\\
        &\cong\bigl(\Pin_n^-\grotimes\sAut_A(M)\bigr)_{\ev}\,,
        \end{aligned}
    \end{equation}
    with its natural projection map  to $O_n$.
    $BH(A;M)\to BO$ is a stable tangential structure, i.e.\  a $H_n(A;M)$-structure on $V$ is equivalent to a $H_{n+1}(A;M)$-structure on $V\oplus \underline{\R}$.
    The tangential structure is oriented, i.e.\ factors as $BH(A;M)\to BSO$, if and only if $M$ is oriented, i.e.\ if $M\ncong \Pi M$.
\end{theorem}
In particular, $(A,M)$-structures can be glued and are ``homotopical''.
This can also be proven directly.

\begin{proof}
    Let $(\R^n, \Sigma_n)$ be the local model of a type $(A,M)$-structure on $\R^n$ (Example~\ref{ex:localmodel}), where $\Sigma_n=\Cl_n\grotimes M$.
    There is an explicit group homomorphism to the group of grading-preserving and grading-reversing automorphisms:
    \begin{equation*}
    \begin{aligned}
        \Pin_n^-\grotimes \sAut_A(M) &\longrightarrow \sAut(\R^n,\Sigma_n)\\
        a\grotimes \phi & \longmapsto \left[\mathrm{Ad}_a, a\grotimes \phi: \Sigma_n \to (\Sigma_n)_{\mathrm{Ad}_a}\right]
    \end{aligned}
    \end{equation*}
    Here $\Pin_n^-\subset(\Cl_{-n})^\times$ acts on $\Cl_n$ from the left by $a\cdot x=(-1)^{|a||x|}xa$, i.e.\ by right multiplication, $\sAut_A(M)$ acts on $M$ by evaluation, and
    \begin{equation*}
        (a\grotimes\phi)\cdot(x\grotimes m)=(-1)^{|\phi||x|}\,a\cdot x\grotimes \phi(m) \,.
    \end{equation*}
    This commutes with the left $A$-action, and it intertwines the right $\Cl(\R^n)$-action with its twist by $\mathrm{Ad}_a(v)=(-1)^{|a|}ava^{-1}$, the twisted adjoint action.
    To show that this is an isomorphism of fermionic groups, we just note that the twisted adjoint action $\mathrm{Ad}:\Pin_n^-\to O(n)$ is surjective with kernel $\{\pm 1\}$ and that $\sAut_{A-\Cl_n}(\Sigma_n)\cong \sAut_A(M)$.
    This implies that for any $(f\in O_n, \Phi:\Sigma_n\to (\Sigma_n)_f)$ there is an $a\in \Pin_n^-$ with $\mathrm{Ad}_a=f$ and a $\phi\in \sAut_A(M)$ such that $\Phi=a\grotimes \phi$:
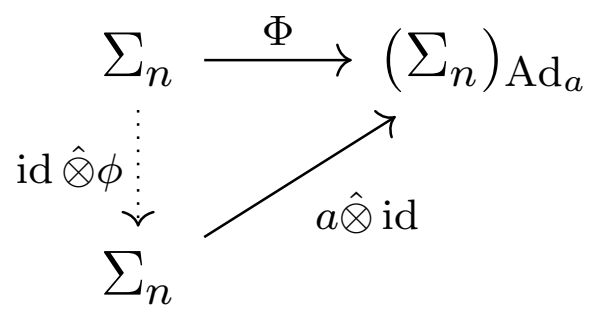
    \begin{equation*}
        \begin{tikzcd}
            \Sigma_n \ar[r,"\Phi"] \ar[d,dotted,"\id\grotimes \phi"'] & (\Sigma_n)_{\mathrm{Ad}_a} \\
            \Sigma_n \ar[ru,"a\grotimes \id"']& 
        \end{tikzcd}
    \end{equation*}
    This proves surjectivity of the map above.
    Note that the kernel is precisely $\{1\grotimes \id, -1\grotimes -\id\}$, which is quotiented out in the fermionic tensor product.

    We conclude that there is a homeomorphism on the even subgroups
    \begin{equation*}
        H_n(A;M)\cong(\Pin_n^-\grotimes \sAut_A(M))_{\ev} \cong \Aut(\R^n,\Sigma_n) \,,
    \end{equation*}
    as the map above defines a continuous bijective homomorphism of Lie groups.
    Under this identification the projection to $O_n$ is $a\grotimes\phi\mapsto \mathrm{Ad}_{a}$.
    The first description in \eqref{eq:grouphom} follows by passing to opposites.

    It follows that a rank $n$ $(A,M)$-structure is patched from the trivial structure $(\R^n,\Sigma_n)$ along continuous transition functions with values in $H_n(A;M)$.
    That is, an $(A,M)$-structure is equivalent to an $H_n(A;M)$-principal bundle and $BH_n(A;M)$ is indeed a classifying space.

    If $M\ncong \Pi M$ is oriented, then there is no bimodule isomorphism $\Sigma_n \to  (\Sigma_n)_f$ for $\det(f)=-1$, since the target is isomorphic to $\Pi \Sigma_n$.
    Hence the projection $\Aut(\R^n,\Sigma_n)\to O(n)$ factors through $SO(n)$ and the tangential structure is oriented.
    The converse is similar.

    We are left to prove stability of the tangential structure.
    Suppose that $V\oplus \underline{\R}^k$ carries an $(A,M)$-structure with $\Sigma$ a $A$-$\Cl(V\oplus \underline{\R}^k\!-\!\rk V\!-\!k)$-bimodule bundle of type $M$.
    Using the Morita equivalence $\Cl(\underline{\R}^k\!-\!k)\simeq \R$, there also exists a $A$-$\Cl(V\!-\!\rk V)$-bimodule bundle $\Sigma'$ of type $M$. 

\end{proof}

Note that \eqref{eq:grouphom} defines a tangential structure $BH(A;M)\to BO$ for every super Banach algebra $A$ and nonzero finitely generated projective module $M$.
For $M=0$ we will set $H_n(A;M)=O_n$. 
\begin{remark}
    \label{rmk:ebert1}
    Ebert~\cite[Ex.~2.11]{EbertIndex} introduces \emph{unstable} tangential structures $\theta_A(d)\to BO(d)$ for graded Real $C^*$-algebras $A$.
    A $\theta_A(d)$-structure on a rank $d$ inner product bundle $V$ is a bundle $Q$ of finitely generated projective graded right Hilbert $A$-modules with an odd $A$-linear map $c:V\to\End_A(Q)$ satisfying $c(v)^2=-|v|^2$.
    In our conventions, $Q$ is thus a $\Cl(-V)$-$A$-bimodule bundle, or equivalently, using $\Cl(-V)^{\op}\cong\Cl(V)$, a bimodule bundle $Q:\Cl(V)\squigrightarrow A^{\op}$.
    Since $\Cl_{-d}$ is Morita invertible, $Q\mapsto Q\grotimes\Cl_{-d}$ is an equivalence, natural in isometries of $V$, onto bimodule bundles $\Cl(V\!-\!d)\squigrightarrow A^{\op}\grotimes\Cl_{-d}$.
    We have
    \[
        \theta_A(d)\simeq\coprod_{[N]} BH_d(A^{\op}\grotimes\Cl_{-d};N),
    \]
    where $[N]$ runs over types of finitely generated projective $A^{\op}\grotimes\Cl_{-d}$-modules.
\end{remark}
\subsection{Equivalence and irreducible bibundles}
\label{sec:equivalencesandirreps}
In this subsection we explain important instances of type $(A,M)$-structures.
We recover the spin and pin groups and the spacetime symmetry groups of
Freed--Hopkins~\cite[Tables~(9.24)--(9.25)]{FreedHopkins}.

\begin{enumerate}
    \item Type $(A,M)$ for $M$ irreducible.
    A bimodule bundle $\Sigma$ of type $(A,M)$ is fiberwise irreducible if and only if $M$ is irreducible.
    By Schur's lemma and Remark~\ref{rmk:handedness},
    \[
        D=\mathrm{sEnd}_A(M)^{\op}
    \]
    is a super division algebra and $M$ is an $A$-$D$-bimodule.
    The functor $M\grotimes_D-$ induces an equivalence between $(D,D)$-structures
    and $(A,M)$-structures; its inverse is
    $\mathrm{sHom}_A(M,-)$.
    This reduces to type $(D,D)$.
   
    \item Type $(A,A)$.
    Every algebra $A$ acts on itself and we call $A$ the left regular $A$ module.
    Theorem~\ref{thm:typedstructures} gives
    \[
        H_n(A;A)
        =(\Pin_n^+\grotimes A^\times)_{\ev}.
    \]
    \item Type $(\R,\R)$.
    A type $(\R,\R)$ bimodule induces a fiberwise Morita equivalence
    $\Cl(V\!-\!\rk V)\simeq\R$.
    Moving the constant Clifford factor across this equivalence gives:
    \begin{proposition}
        \label{prop:spinstructures}
        For a real vector bundle $V$, the following data
        are equivalent:
        \begin{enumerate}
            \item a $\Spin$-structure on $V$;
            \item a type $(\R,\R)$ bimodule bundle;
            \item a Morita equivalence between $\Cl(V)$ and
            the constant algebra bundle $\Cl_{\rk V}$;
            \item a fiberwise irreducible $\Cl(V)$-$\Cl_{\rk V}$-bimodule bundle.
        \end{enumerate}
    \end{proposition}
    Here and below the noncompact automorphism groups can be
    replaced by their compact reductions without changing the
    homotopy type of the classifying space.

    \item 
    Type $(D,D)$ for a super division algebra $D$.
    There are ten real super division algebras and they all admit
    compatible graded $C^*$-structure. The compact homogeneous
    unit group is
    \[
        O(D)=\{d\in D_0\cup D_1:d^*d=1\}.
    \]
    Polar decomposition gives a homotopy equivalence $D^\times\simeq O(D)$.
    This is the union of the unit spheres in the homogeneous
    summands.

    For $s\in\Z/8$, let $D(s)$ be the real super division
    algebra in the Morita class of $\Cl_s$.
    The two complex super division algebras are $D^c(0)=\C$ and
    $D^c(1)=\Cxl_{-1}$, following~\cite[(9.24)]{FreedHopkins}.
    \begin{definition}
        \label{def:fermionictypes}
        The spacetime groups are
        \[
            H_n(s)=(\Pin_n^+\grotimes O(D(s)))_{\ev},
            \qquad
            H_n^c(s)=(\Pin_n^+\grotimes O(D^c(s)))_{\ev}.
        \]
    \end{definition}
    These are compact models for $H_n(D(s);D(s))$ and its complex
    analogue. They have natural maps to $O_n$ and to $O(D(s))_b$,
    induced by $(K\grotimes L)_b=K_b\times L_b$.
    With these conventions the ten groups are
    \begin{equation}
    \label{eq:divisiongroups}
    \begin{array}{c l c c}
    \toprule
    D&H&O(D)_b&s\\ \midrule
    \R&\Spin&*&0\\
    \Cl_{-1}&\Pin^+&O_1&-1\\
    \Cl_{-2}&\Pin^{\tilde c+}&O_2&-2\\
    \Cl_{-3}&\Pin^{h-}&O_3&-3\\
    \mathbb H&\Spin^h&SO_3&4\\
    \Cl_{+3}&\Pin^{h+}&O_3&3\\
    \Cl_{+2}&\Pin^{\tilde c-}&O_2&2\\
    \Cl_{+1}&\Pin^-&O_1&1\\
    \C&\Spin^c&U_1&0^c\\
    \Cxl_{-1}&\Pin^c&U_1\times O_1&1^c\\
    \bottomrule
    \end{array}
    \end{equation}
    We specify auxiliary virtual bundles needed below. Put $BO(D)_b=B(O(D)/\{\pm1\})$.
    Let $L_k\to BO_k$, $Q\to BSO_3$, $E\to BU_1$ and $\lambda\to BO_1$ be the universal bundles.
    Define rank-zero virtual real vector bundles over $BO(D)_b$ by
    \begin{equation}
        \label{eq:auxiliarybundles}
        \xi_D=
        \begin{cases}
            0,&D=\R,\\
            L_k\ominus k,&D=\Cl_{+k},\quad k=1,2,3,\\
            k\ominus L_k,&D=\Cl_{-k},\quad k=1,2,3,\\
            Q\ominus 3,&D=\mathbb H,\\
            E_{\R}\ominus 2,&D=\C,\\
            E_{\R}\oplus\lambda\ominus 3,&D=\Cxl_{-1}.
        \end{cases}
    \end{equation}

    \begin{lemma}[{\cite[(10.18)--(10.22) and (10.36)]{FreedHopkins}}]
        \label{lem:pullbackfermionicgrps}
        There is a homotopy pullback
        \[
        \begin{tikzcd}
            BH_n(D)\ar[r]\ar[d]&B\Spin\ar[d]\\
            BO_n\times BO(D)_b\ar[r,"-\oplus\xi_D"']&BO.
        \end{tikzcd}
        \]
        Here $H_n(D)=(\Pin_n^+\grotimes O(D))_{\ev}$.
        Consequently there is a stable splitting
        $BH(D)\simeq B\Spin\times BO(D)_b$.
    \end{lemma}

    \begin{proposition}
        \label{prop:typeddbibundles}
        For a rank $n$ vector bundle $V$, the following data are
        equivalent, after compact reduction:
        \begin{enumerate}[label=(\alph*), ref=(\alph*)]
            \item an $H_n(D)$-structure on $V$;
            \item a type $(D,D)$ bimodule bundle;
            \item a fiberwise irreducible
            $\Cl(V)$-$(\Cl_n\grotimes D)$-bimodule bundle;
            \item \label{item:auxiliarybundle} a principal $O(D)_b$-bundle $P$ and a spin
            structure on $V\oplus \xi_D(P)$.
        \end{enumerate}
    \end{proposition}
    \begin{proof}
        The first two descriptions follow from
        Theorem~\ref{thm:typedstructures}.
        Dualizing a $D$-$\Cl(V\!-\!n)$-bimodule and moving the
        constant $\Cl_{-n}$ factor to the right gives the third
        description. These operations preserve irreducibility;
        over a super division algebra an irreducible module is
        the regular module up to parity.
        The final description is Lemma~\ref{lem:pullbackfermionicgrps}.
    \end{proof}
    \begin{example}
        For $D=\Cl_{\mp1}$, Proposition~\ref{prop:typeddbibundles}\ref{item:auxiliarybundle} is the usual description of a
        $\Pin^\pm$-structure as a
        $\Spin$-structure on $V\mp\det V$~\cite[Lem.~1.7]{KirbyTaylor}.
        Indeed, the degree condition identifies the auxiliary
        line with $\det V$.
        The obstruction is $w_2(V)$ for $\Pin^+$ and
        $w_2(V)+w_1(V)^2$ for $\Pin^-$~\cite[Lem.~1.3]{KirbyTaylor}.
    \end{example}
    \begin{example}
        A $\Spin^c$-structure consists of a complex line $L$ and
        a spin structure on $V\oplus L_{\R}$.
        Its existence is equivalent to
        \[
            w_1(V)=0,\qquad
            w_2(V)=c_1(L)\bmod2.
        \]
        For oriented $V$, such a line exists exactly when the third integral Stiefel--Whitney class
        $W_3(V)=\beta(w_2(V))$ vanishes.
    \end{example}
    \begin{proof}[Proof of Lemma~\ref{lem:pullbackfermionicgrps}]
    The fermionic group structure provides a map $\kappa_D:BO(D(s))_b\to \tau_{\leq 2}BO$.
    By example~\ref{ex:pinfermionic}, the fermionic group structure on $\Pin_n^+$ is the truncation map $BO_n\to \tau_{\leq 2}BO$.
    One checks that in all cases the following diagram commutes:
    \begin{equation}
    \label{eq:1}
        \begin{tikzcd}
            BO_n\times BO(D(s))_b \ar[r,"{(\tau,\kappa_D)}"]\ar[d,"V\oplus \xi_D"] & \tau_{\leq 2}BO \times \tau_{\leq 2} BO \ar[d,"\oplus"] \\
            BO \ar[r] &\tau_{\leq 2}BO
        \end{tikzcd}
    \end{equation}
    By definition, the diagonal composite in square~\eqref{eq:1} classifies the fermionic group $\Pin_n^+\grotimes O(D(s))$.
    Its even part is the fiber of this composite.
    By pullback pasting, the fiber of a composite is the pullback of the fiber of the second map.
    We arrive at the claim by noting that $\fib(BO\to \tau_{\leq 2}BO)=B\Spin$.
\end{proof}
\end{enumerate}

\subsection{Existence of irreducible Clifford bundles}
\label{sec:existenceirreducible}
The existence of an irreducible Clifford module bundle does not in
general imply the existence of a spin structure.
Chen~\cite{xuanchen} exhibits counterexamples to this for oriented Riemannian
manifolds using the negative-definite algebra $\Cl(-TX)$.

\begin{center}
    \emph{When does an irreducible Clifford module bundle exist?}
\end{center}
We first give the positive-definite statement for general vector bundles $V$.
It is possible to specialize to $V=TX$.
\begin{proposition}
    \label{prop:irreduciblecliffordbundles}
    Let $V$ be a real vector bundle of rank $n$ with inner product.
    \begin{enumerate}
        \item A graded irreducible $\Cl(V)$-module bundle exists
        if and only if $V$ admits an $H(-n)$-structure.
        \item An ungraded irreducible $\lvert\Cl(V)\rvert$-module bundle exists
        if and only if $V$ admits an $H(1-n)$-structure.
    \end{enumerate}
    The indices are read modulo $8$; the groups $H(s)$ are listed in~\eqref{eq:divisiongroups}.
    For complex Clifford modules the same statements hold with
    $H^c$ and indices modulo $2$.
\end{proposition}
\begin{proof}
    A graded $\Cl(V)$-module is a $\Cl(V)$-$\R$-bimodule.
    Since $\Cl_n\grotimes D(-n)$ is Morita equivalent to $\R$,
    the first assertion follows from
    Proposition~\ref{prop:typeddbibundles}.

    For the second assertion, use the equivalence between ungraded and graded modules
    \[
        \Mod_{\lvert\Cl(V)\rvert}
        \simeq \Mod^{C_2}_{\Cl(V)\grotimes\Cl_{-1}}.
    \]
    Explicitly, an ungraded module $W$ with action $\rho$ is sent
    to $W\oplus W$, with the two summands as its grading and actions
    \[
        v\longmapsto
        \begin{pmatrix}0&\rho(v)\\ \rho(v)&0\end{pmatrix},
        \qquad
        f\longmapsto
        \begin{pmatrix}0&-1\\1&0\end{pmatrix}.
    \]
    Conversely the invertible odd operator $f$ identifies the
    two graded summands and recovers the ungraded module.
    
    Since $\Cl_{-1}=\Cl_{+1}^{\op}$, the resulting module
    is a $\Cl(V)$-$\Cl_{+1}$-bimodule.
    The coefficient algebra $\Cl_n\grotimes D(s)$ has the required
    Morita class precisely when $n+s=1$ in $\Z/8$.
    The complex argument is identical.
\end{proof}

\begin{remark}[Negative-definite convention]
    \label{rmk:negativeclifford}
    The existence of a graded irreducible $\Cl(-V)$-module is equivalent to the existence of a 
    $H(-n)$-structure, whereas the criterion for an ungraded irreducible module is a lift to $H(-n-1)$.
\end{remark}
\begin{remark}
    For $n=0,6,7\bmod8$, the negative-definite ungraded types are
    respectively $-1,1,0$, corresponding to $\Pin^+$, $\Pin^-$ and
    $\Spin$. On an oriented bundle either $\Pin$-structure restricts to
    a $\Spin$-structure. This recovers~\cite[Theorem~4]{xuanchen}:
    In these dimensions, granted that $V$ is oriented, the existence of an ungraded irreducible $\Cl(-V)$-module bundle is equivalent to the existence of a spin structure on $V$.
    Over $\C$, each oriented restriction of $H^c(0)$ or $H^c(1)$
    is $\Spin^c$. Thus, on an oriented manifold, the existence of
    either a graded or an ungraded irreducible complex Clifford module
    bundle is equivalent to the existence of a $\Spin^c$-structure,
    recovering~\cite[Theorem~2]{xuanchen} in the ungraded case.
    
\end{remark}

\subsection{The Picard groupoid of super Banach algebras}
\label{sec:picardgroupoid}
Every finite dimensional super algebra defines a super Banach algebra by equipping it with a norm.
We will now briefly recall the setup of the Morita $(\infty,1)$-category $\Bim(\sBan_\R)_\infty$ of real super Banach algebras from \cite{aretzstehouwer}.
\begin{definition}
    Let $\Bim(\sBan_\R)$ be the symmetric monoidal $(2,1)$-category of real super Banach algebras, bimodules which are finitely generated projective as a module over the target algebra and invertible intertwiners.
    The monoidal structure is the graded projective tensor product $\grotimes$.
    A bimodule homotopy is $M:A\squigrightarrow C([0,1];B)$ regarded as a homotopy of bimodules $\ev_i\circ M:A\squigrightarrow B$.
    A bimodule homotopy equivalence is a bimodule $M:A\squigrightarrow B$ such that there exists a bimodule $N:B\squigrightarrow A$ and bimodule homotopies $N\grotimes_B M \simeq A$ and $M\grotimes_A N \simeq B$.
    We let $\Bim(\sBan_\R)_\infty=\Bim(\sBan_\R)[W^{-1}]$ be the $(\infty,1)$-category obtained by localizing at bimodule homotopy equivalences.
\end{definition}
The mapping spaces in $\Bim(\sBan_\R)_\infty$ should be interpreted as \emph{spaces of bimodules}. 
By the graded version of \cite[Thm.~2.101]{aretzstehouwer} we have
\begin{equation}
    \label{eq:mappingspacesbimban1}
            \Hom_{\Bim(\sBan_\R)_\infty}(A,B) \simeq \bigsqcup_{[P]} \Alg(A,\sEnd_B(P)^{\op}) \sslash \mathsf{Aut}_B(P)
            \,,
\end{equation}
where $[P]$ runs over isomorphism classes of finitely generated projective $B$-modules and where $\sslash$ denotes the homotopy quotient by the conjugation action on the topological space of grading-preserving Banach algebra homomorphisms equipped with the $k$-ification of the compact-open topology~\cite[Def.~3.10]{aretzstehouwer}.
If $A$ is finite-dimensional, then this agrees with the norm topology.

\begin{remark}
    Every grading-preserving algebra homomorphism induces a bimodule and this assembles into a symmetric monoidal functor 
    \[
        \Alg(\sBan_\R)_\infty \to \Bim(\sBan_\R)_\infty 
    \] 
    from the $\infty$-category $\Alg(\sBan_\R)_\infty=\Alg(\sBan_{\R})[W^{-1}]$ of super Banach algebras.
    The mapping spaces are spaces of algebra homomorphisms topologized as above.
\end{remark}

\subsection{A categorical Postnikov tower}
\label{sec:categoricalpostnikov}
We will now show that the Clifford algebras categorify the determinant in a precise sense.
While the determinant can basically be identified with the first Stiefel--Whitney class, the Clifford algebras are a categorification of the first and second Stiefel--Whitney classes.

\begin{figure}[htbp]
    \centering
    \begin{tikzpicture}[>=stealth, every node/.style={inner sep=3pt}]
        \node (vspin) at (0,1.6) {$\Vect^{\Spin}$};
        \node (pic) at (3.6,1.6) {$\Pic(\Bim(\sBan_\R)^{\mathrm{fd}}_\infty)$};
        \node (vso) at (0,0) {$\Vect^{SO}$};
        \node (vo) at (1.8,0) {$\Vect^O$};
        \node (sline) at (3.6,0) {$\sLine_{\R}$};
        \draw[->] (vspin) -- (vso);
        \draw[->] (vspin) -- (vo);
        \draw[->] (vso) -- (vo);
        \draw[->] (vo) -- node[below] {$\det$} (sline);
        \draw[->] (vo) -- node[midway, sloped, above] {$\Cl$} (pic);
        \draw[->,dotted] (pic) -- node[right] {$\operatorname{trace}$} (sline);
    \end{tikzpicture}%
    \hspace{1cm}%
    \begin{tikzpicture}[>=stealth, every node/.style={inner sep=3pt}]
        \node (bspin) at (0,1.6) {$B\Spin$};
        \node (postnikov) at (3.6,1.6) {$\tau_{\leq 2}BO$};
        \node (bso) at (0,0) {$BSO$};
        \node (bo) at (1.8,0) {$BO$};
        \node (bz2) at (3.6,0) {$B\Z/2$};
        \draw[->] (bspin) -- (bso);
        \draw[->] (bspin) -- (bo);
        \draw[->] (bso) -- (bo);
        \draw[->] (bo) -- node[below] {$w_1$} (bz2);
        \draw[->] (bo) -- node[midway, sloped, above] {$(w_1,w_2)$}(postnikov);
        \draw[->] (postnikov) -- (bz2);
    \end{tikzpicture}
    \caption{The categorical Postnikov/Whitehead tower (left) and its stable homotopical
    counterpart (right), with the corresponding Postnikov stages in the right-hand column of each.}
    \label{fig:categorical-whitehead-tower-postnikov}
\end{figure}
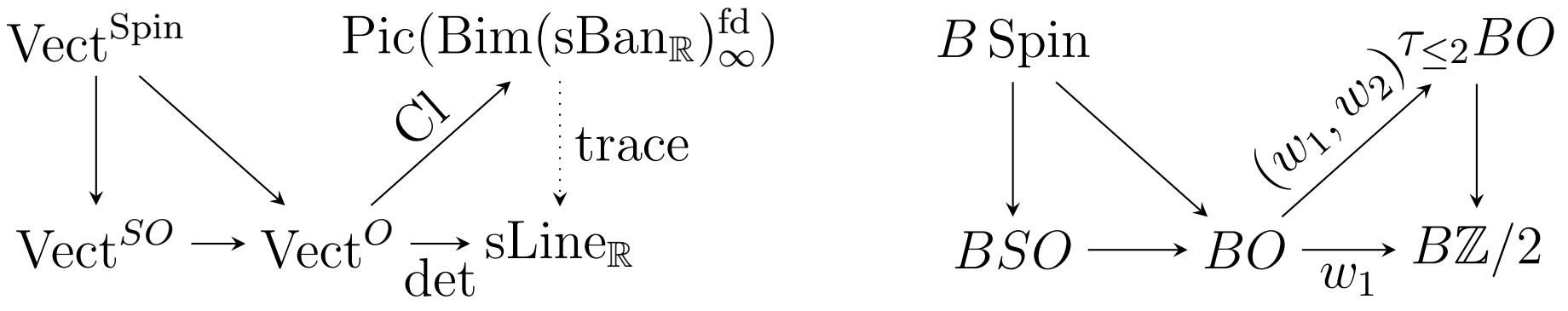
The Clifford algebra construction $V\mapsto \Cl(V)$ defines a symmetric monoidal functor~\cite[Prop.~3.39]{aretzstehouwer}\footnote{Here, we use the \emph{positive-definite} convention, where $v^2=\langle v,v \rangle$.}
from inner product spaces to the Morita $(\infty,1)$-category
\begin{equation}
    \Cl:\Inn \to \Alg(\sBan_{\R})_\infty \to \Bim(\sBan_\R)_\infty \,.
\end{equation}
which restricts to $\Cl:\Vect_{\R}^{\cong}\simeq \bigsqcup BO_n \to \Pic(\Bim(\sBan_\R)^{\mathrm{fd}}_\infty)$ taking values in the Picard $\infty$-groupoid of \emph{finite-dimensional semisimple} $\grotimes$-invertible algebras and invertible bimodules.
Finite-dimensional semisimple super algebras are fully dualizable with dual $A^{\op}$.

We also recall the description~\cite[Thm.~3.65]{aretzstehouwer}: The Picard $\infty$-groupoid
$\Pic(\Bim(\sBan_\R)^{\mathrm{fd}}_\infty)$ is a $2$-type that has the following description:
        The equivalence classes of objects are generated by the Clifford algebras:
        \begin{equation*}
            \pi_0\Pic(\Bim(\sBan_\R)^{\mathrm{fd}}_\infty)=\{\Cl_0,\dots,\Cl_7\}\cong \Z/8 \,.
        \end{equation*}
        The loop space $\Omega_{\mathbbm{1}}(\Pic(\Bim(\sBan_\R)^{\mathrm{fd}}_\infty))$ is the subspace of $\Hom_{\Bim(\sBan_\R)_\infty}(\R,\R)$ on the invertible bimodules $\R$ and $\Pi \R$.
        This is the   $1$-type $\sLine_{\R}$ of real graded lines with
        \begin{align*}
            \pi_1\Pic(\Bim(\sBan_\R)^{\mathrm{fd}}_\infty)&=\pi_0 N^{\mathrm{hc}}(\sLine_{\R}^{\cong})=\{\R, \Pi \R\} \cong \Z/2 \\
            \pi_2\Pic(\Bim(\sBan_\R)^{\mathrm{fd}}_\infty)&=\pi_1 N^{\mathrm{hc}}(\sLine_{\R}^{\cong})=\pi_0\frac{\{\text{linear isom. } \R\to \R\}}{\text{homotopy}}\cong\Z/2 \,.
        \end{align*}
All higher homotopy groups vanish.
Note that the homotopy groups do not depend on the basepoint, since the Picard $\infty$-groupoid is grouplike.

We can now prove that $\Cl$ is a sort of Postnikov truncation of $BO$.
\begin{proposition} 
    \label{prop:cliffordpostnikov}
    The symmetric monoidal functor $\Cl$ restricts to a map of $\E_\infty$-spaces
    \begin{equation*}
        \Cl:\Inn^{\cong}\simeq \bigsqcup BO_n \to \Pic(\Bim(\sBan_\R)^{\mathrm{fd}}_\infty) \,.
    \end{equation*}
    On $\pi_0$ this map is given by the mod $8$ dimension of the vector space and induces the surjection $\N \to \Z/8$. It is an isomorphism on $\pi_1$ for all components $n\geq 1$ and on $\pi_2$ for all $n\geq 3$. For $n=2$, the map on $\pi_2$ is reduction modulo $2$, $\Z\twoheadrightarrow\Z/2$.
    The target is grouplike and hence the map extends to the group completion $\Inn^{\cong,\gp} \simeq \Z\times BO$. 
\end{proposition}
\begin{proof}
    The map clearly restricts.
    From the description of the Picard $\infty$-groupoid above, we are only left to prove that it induces isomorphisms on $\pi_1,\pi_2$.
    For $n\geq 1$, the two components of $\Omega BO_n=O_n$ are distinguished by the determinant. Example~\ref{ex:irreducibleClnnmodules} identifies the $\Cl_n$-bimodule $(\Cl_n)_f$ associated to $f\in O_n$ with $\Pi^{\epsilon(f)}\Cl_n$, where $\epsilon(f)=(1-\det(f))/2\in\{0,1\}$. This proves the assertion on $\pi_1$.

    Put $A=\Cl(\R^n)$.
    Let $\gamma_t$ be a path in $SO_n$ starting at the identity, and let $\widetilde{\gamma}_t$ be the lift to $\Spin_n\subset \Cl(\R^n)^{\ev}$ with $\widetilde{\gamma}_0=1$.
    The induced algebra automorphism is $\Cl(\gamma_t)=\mathrm{Ad}_{\widetilde{\gamma}_t}$, and the associated family of bimodules is
    \begin{equation*}
        M_t=A_{\Cl(\gamma_t)},\qquad a\cdot m\cdot_t b=a\,m\,\widetilde{\gamma}_t b \widetilde{\gamma}_t^{-1}.
    \end{equation*}
    Equivalently, its sections form $C([0,1];A)$ with this pointwise action.
    The spin lift gives an explicit continuous trivialization by even bimodule isomorphisms
    \begin{equation*}
        T_t:M_t\xrightarrow{\cong}A,\qquad T_t(m)=m\,\widetilde{\gamma}_t,
    \end{equation*}
    since
    \begin{equation*}
        T_t(a\,m\cdot_t b)
        =a\,m\,(\widetilde{\gamma}_t b \widetilde{\gamma}_t^{-1})\widetilde{\gamma}_t
        =a\,T_t(m)\,b.
    \end{equation*}
    The nontrivial loop in $\pi_1SO_n$ for $n\geq 3$ lifts to a path in $\Spin_n$ with $\widetilde{\gamma}_1=-1$, i.e.\ the nontrivial bimodule isomorphism $A\cong A$, which generates $\pi_2\Pic(\Bim(\sBan_\R)^{\mathrm{fd}}_\infty)$ (of the path component of $A=\Cl(\R^n)$).    
\end{proof}
We also write $\Cl$ for the stable Clifford functor $\Cl:BO \to \Pic(\Bim(\sBan_\R)^{\mathrm{fd}}_\infty)$.
\begin{corollary}
    \label{cor:spinfiberseq}
    There are fiber sequences of $\E_\infty$-spaces
    \begin{equation}
        B\Spin \to BO \xrightarrow{\Cl} \Pic(\Bim(\sBan_\R)^{\mathrm{fd}}_\infty) \,.
    \end{equation}
    \begin{equation}
        B\Spin^c \to BO \xrightarrow{\Cxl} \Pic(\Bim(\sBan_\C)^{\mathrm{fd}}_\infty) \,.
    \end{equation}
\end{corollary}
\begin{proof}
    The first fiber sequence follows directly by the definition of $B\Spin$ as a $2$-connected cover of $BO$.
    The second one uses a similar description of the complex Clifford functor $\Cxl$ as $(w_1,W_3)$ into the Picard groupoid of complex super algebras. 
\end{proof}
This explains the definition (Proposition~\ref{prop:spinstructures}) of a $\Spin$-structure on $V$ as a Morita equivalence $\Cl(V)\simeq \Cl(\rk V)$ through the lens of homotopy theory: A lift to $B\Spin$ is precisely a trivialization of $\Cl(V)$ in the Morita $(\infty,1)$-category $\Bim(\sBan_\R)_\infty$.

\section{The natural transformation \texorpdfstring{$J\to K^{\gr}\circ \Cl$}{J->K(Cl)}}
\label{sec:naturaltransformation}
Let $\Inn$ be the symmetric monoidal $\infty$-category of finite-dimensional inner product spaces and isometric embeddings.
The monoidal structure is direct sum $\oplus$.
Building on the construction of graded $K$-theory from \cite{aretzstehouwer},
we will construct a monoidal natural transformation of lax symmetric monoidal functors $\Inn\to \Sp_{\geq 0}$
\begin{equation}
\label{eq:monoidaltransformation}
        S^V\longrightarrow k^{ABS}_{\Cl(V)}\,,
\end{equation}
which we compose with the map $k^{ABS}\to K^{\gr}=k^{ABS}[\beta^{-1}]$ to obtain
\begin{equation*}
    J(V)=S^V \longrightarrow K^{\gr}_{\Cl(V)} \,.
\end{equation*}
We then show that this induces an equivalence of $KO$-modules
\begin{equation*}
    S^V\otimes KO \xrightarrow{\simeq} K^{\gr}_{\Cl(V)} \,.
\end{equation*}
In particular, $K^{\gr}:\Pic(\Bim(\sBan_\R)^{\mathrm{fd}}_\infty)\to \Pic(\Mod(KO))$ is symmetric monoidal.
This will prove Theorem~\ref{thm:B}.

Recall the definition of connective ABS $K$-theory $k^{ABS}$ from \cite[Def.~3.42]{aretzstehouwer}:
For any super Banach algebra $A$, we let
\begin{equation*}
    k_A^{ABS}= \cofib \left(k^{C_2}_{A\grotimes \Cl_{-1}}\xrightarrow{i_-^*} k^{C_2}_A  \right) \,,
\end{equation*}
where 
\begin{equation*}
    k^{C_2}_A=(\Mod^{C_2}_A)^{\cong,\gp}
\end{equation*}
is the group completion of the $\E_\infty$-monoid of finitely generated projective left $A$-modules and $i_-^*$ is induced by forgetting the $\Cl_{-1}$-action.

\subsection{The \texorpdfstring{$J$}{J}-homomorphism}
Recall from \cite[Cor.~3.30]{aretzstehouwer} that the arrow $\R^1\xrightarrow{\op} 0$ in $\Inn^{\op}$ refines to a Smith ideal, i.e.\ a lax symmetric monoidal functor from the arrow category $[1]_{\mathrm{min}}$ with the minimum symmetric monoidal structure.
The (co)Yoneda embedding
\begin{align*}
    \Inn^{\op}&\longrightarrow \Fun(\Inn,\Spc)\\
    V& \longmapsto [V^*:W\mapsto \Inn(V,W)]
\end{align*}
is symmetric monoidal with respect to the Day convolution symmetric monoidal structure on the target~\cite[Cor.~4.8.1.12, Rem.~4.8.1.13]{HA}, see also the proof of Lemma~\ref{lem:mappingspaces}.
The cofiber/quotient of a Smith ideal is a commutative algebra object.
Commutative algebra objects with respect to Day convolution are lax symmetric monoidal functors.
\begin{definition}
    The $J$-homomorphism is the (lax) symmetric monoidal functor $\Inn\to \Spc_*$ defined by
    \begin{equation}
    J:=\cofib(\R^*\to 0^*) \in \CAlg(\Fun(\Inn,\Spc_*)) \,.
    \end{equation}
\end{definition}
The underlying functor of $J$ is 
\begin{equation*}
    J(V)=\cofib\left(\Inn(\R,V)\to \Inn(0,V)\right) = \cofib(S(V)\to *)=S^V \,.
\end{equation*}
$J(V)$ comes with a natural $O(V)$-action which readily identifies with the natural action on the $1$-point compactification of $V$.
Since $J(V)\wedge J(W)\simeq J(V\oplus W)$, $J$ is strong symmetric monoidal.
\begin{remark}
    We will also write $J$ for the symmetric monoidal composite 
    \begin{equation*}
        \Inn \xrightarrow{J} \Spc_* \xrightarrow{\Sigma^\infty} \Sp \,.
    \end{equation*}
\end{remark}

\subsection{Construction}
We start with an easy technical lemma whose proof we defer to Appendix~\ref{app:mappingspaces}.
\begin{lemma}[Lemma~\ref{lem:laxtransformationapp}]
    \label{lem:laxtransformation}
    Let $\calC$ and $\calD$ be symmetric monoidal $\infty$-categories such that every object $d\in \calD$ is dualizable with dual $d^\vee$.
    Let $F:\calC\to \calD$ be symmetric monoidal and let $G:\calD\to \Sp$ be lax symmetric monoidal.
    Then there is a lax symmetric monoidal natural transformation of functors $\calC^{\op}\times \calC\to \Sp$
    \begin{equation*}
        \Sigma_+^\infty \Map_{\calC}(c,c')\longrightarrow G(F(c)^{\vee}\otimes F(c'))\,.
    \end{equation*}
    Its adjoint $\Map_{\calC}(c,c')\to \Omega^\infty G(F(c)^{\vee}\otimes F(c'))$ sends $f:c\to c'$ to the image under $G$ of the morphism $\mathbbm{1}\to F(c)^{\vee}\otimes F(c')$ corresponding to $F(f)$, precomposed with the unit $\mathbb{S}\to G(\mathbbm{1})$.
\end{lemma}
We consider $\calC=\Inn$, $\calD=\Bim(\sBan_\R)^{\mathrm{fd}}_\infty$, $F=\Cl$ and $G=k^{C_2}$ to obtain a lax symmetric monoidal natural transformation of functors $\Inn^{\op}\times \Inn \to \Sp$:
\begin{equation}
    \label{eq:naturaltransformation}
    \Sigma_+^\infty \Inn(V,W) \to k^{C_2}_{\Cl(-V)\grotimes \Cl(W)}
\end{equation}
Restricting along the Smith ideal $[1]_{\mathrm{min}}\to \Inn^{\op}$ we obtain a lax symmetric monoidal natural transformation of functors $[1]_{\mathrm{min}}\times \Inn\to \Sp$.
At $V\in \Inn$ this is the following commutative square:
\begin{equation}
    \label{eq:deflaxmonoidaltrans}
    \begin{tikzcd}
        \Sigma_+^\infty \Inn(\R,V) \ar[d]\ar[r]& k^{C_2}_{\Cl(V)\grotimes \Cl_{-1}}\ar[d,"i_-^*"]\\
        \Sigma_+^\infty \Inn(0,V) \ar[r]& k^{C_2}_{\Cl(V)}
    \end{tikzcd}
\end{equation}
Since taking cofibers is symmetric monoidal we obtain a lax symmetric monoidal natural transformation of functors $\Inn\to \Sp$:
\begin{equation*}
    J(V) \longrightarrow k^{ABS}_{\Cl(V)} \,.
\end{equation*}

\begin{remark}
    \label{rmk:explicittrans}
    Let us unravel what the adjoint of \eqref{eq:naturaltransformation} $\Inn(V,W) \to \Omega^\infty k^{C_2}_{\Cl(-V)\grotimes \Cl(W)}$ does.
    An isometric embedding $f:V\to W$ induces an algebra homomorphism $\Cl(f):\Cl(V)\to \Cl(W)$ and a further a $\Cl(W)$-$\Cl(V)$-bimodule $\Cl(W)_{\Cl(f)}$.
    By duality this is equivalently a graded left $\Cl(-V)\grotimes \Cl(W)$-module and defines an element of $k^{C_2}_{\Cl(-V)\grotimes \Cl(W)}$.
\end{remark}
\subsection{Proof}
\begin{theorem}
    \label{thm:spin-orientation-proof}
    There is an invertible lax symmetric monoidal natural transformation making the following diagram of lax symmetric monoidal functors commute:
    \begin{equation*}
        \begin{tikzcd}
            \mathcal{I}nn \ar[r,"J"]\ar[d,"\Cl"] &\Sp \ar[d,"\otimes KO"]\\
        \Bim(\sBan_\R)_\infty \ar[r,"K^{\gr}"] & \Mod(KO)
        \end{tikzcd}
    \end{equation*}
    Upon taking cores, it restricts to a symmetric monoidal natural equivalence of symmetric monoidal functors
    \begin{equation*}
        \begin{tikzcd}
            \bigsqcup_n BO_n\simeq \mathcal{I}nn^{\cong} \ar[r,"J"]\ar[d,"\Cl"] &\Pic(\Sp) \ar[d,"\otimes KO"]\\
        \Pic(\Bim(\sBan_\R)^{\mathrm{fd}}_\infty) \ar[r,"K^{\gr}"] & \Pic(\Mod(KO))
        \end{tikzcd}
    \end{equation*}
\end{theorem}
\begin{proof}
    We have constructed a lax symmetric monoidal natural transformation $J\to k^{ABS}\circ \Cl$.
    We can compose with the canonical map $k^{ABS}\to K^{\gr}$ and $KO$-linearize to obtain a lax symmetric monoidal natural transformation $\alpha_V:J(V)\otimes KO \to K^{\gr}(\Cl(V))$.
    Both functors are symmetric monoidal, so it suffices to show that $\alpha$ is an equivalence on $V=\R^1$.
    By Remark~\ref{rmk:explicittrans}, $\alpha$ is induced by
    \begin{equation*}
        \begin{tikzcd}
            S^0=\Inn(\R,\R) \ar[d]\ar[r]& k^{C_2}_{\Cl_{+1}\grotimes \Cl_{-1}} \ar[d]\\
            *= \Inn(0,\R) \ar[r] & k^{C_2}_{\Cl_{+1}}
        \end{tikzcd}   
    \end{equation*}
    and we can explicitly describe the image: The two maps $\pm \id:\R \to \R$ are mapped to the two irreducible $\Cl_{1,1}$-modules, cf.\ Example~\ref{ex:irreducibleClnnmodules}.
    The map $0\to \R$ is mapped to $\Cl_{+1}$ as a module over itself, which is also irreducible.
    Using that $\Mod^{C_2}_{\Cl_{1,1}}$ is equivalent to $\Vect\times \Vect$ and that $\Mod^{C_2}_{\Cl_{+1}}\simeq \Vect$ via the irreducibles above, we see that
    the induced $ko$-linear horizontal maps
    \begin{equation*}
        \begin{tikzcd}
            ko\oplus ko \ar[d]\ar[r]& k^{C_2}_{\Cl_{+1}\grotimes \Cl_{-1}} \ar[d]\\
            ko \ar[r] & k^{C_2}_{\Cl_{+1}}
        \end{tikzcd}
    \end{equation*}
    are equivalences.
    Inverting the Bott element preserves these equivalences.
    Upon taking vertical cofibers we obtain that $\alpha$ induces an equivalence $\Sigma KO \simeq K^{\gr}_{\Cl_{+1}}$.
\end{proof}
\begin{remark}
    $\alpha_{\R^1}$ corresponds to an element $\alpha\in \pi_1K^{\gr}(\Cl_1)\cong \Z$.
 This coincides with the image of the periodicity element $\alpha\in\pi_1k^{ABS}_{\Cl_{+1}}$ constructed in \cite[Sec.~3.3.2, Prop.~3.59]{aretzstehouwer}.
\end{remark}

\section{Applications}

\subsection{Thom spectra}
\label{sec:thomspectra}
Let $\xi:X\to BO$ be the classifying map of a stable vector bundle and let $R$ be a ring spectrum.
We recall from \cite{AndoBlumbergGepner} that $\xi$ is \emph{$R$-oriented} if the composite map 
\begin{equation}
    \label{eq:orientationofvectorbundle}
        X\to BO \xrightarrow{J} \Pic(\Sp) \xrightarrow{\otimes R} \Pic(R)    
\end{equation}
which classifies $R$-homology of the Thom space is homotopic to the constant map $\const_{R}$.
The datum of a nullhomotopy of \eqref{eq:orientationofvectorbundle} is equivalent to a map of spectra
\begin{equation}
    \label{eq:thomspaceorientation}
    M\xi := \colim_X (J\circ \xi) \to R \,,
\end{equation}
which is a \emph{Thom class}, i.e.\ for every $x\in X$
the composite $\mathbb{S}\simeq J(\xi_x)\to M\xi \to R$ is a unit in $\pi_0 R$.
If a universal bundle $\xi:BG\to BO$ admits an $R$-orientation, we say conversely that $R$ is \emph{$G$-oriented}.
In this case we obtain a map
\begin{equation}
    MG := \colim_{BG} J\circ \xi \to R \,.
\end{equation}
We record the $\E_\infty$-version, where rings maps are automatically Thom classes.
\begin{proposition}[{\cite[Thm.~3.5, Lem.~3.15]{AntolinCamarenaBarthel}}]
    \label{prop:orientationshomotopies}
    Let $X$ be a grouplike $\E_\infty$-space and let $\xi:X\to BO$ be an $\E_\infty$-map.
    Then $\E_\infty$-nullhomotopies of \eqref{eq:orientationofvectorbundle} are the same as $\E_\infty$-maps $M\xi\to R$.
\end{proposition}

\begin{remark}
    \label{rmk:mapfromthomspace}
    Let $N$ be an $R$-module spectrum.
    Then a map of spectra $M\xi \to N$ is equivalent to a natural transformation in $\Fun(X,\Mod(R))$ from 
    \[
        X\xrightarrow{\xi} BO \xrightarrow{J} \Pic(\Sp) \xrightarrow{\otimes R} \Mod(R)
    \]    
    to the constant functor $\const_N$.
\end{remark}
\begin{remark}
    For $\xi:X\to BO_n$ a vector bundle, the Thom space $\Th(\xi)$ is also given by the colimit over $X$ of the unstable $J$-homomorphism
    \[
        X\xrightarrow{\xi} BO_n \xrightarrow{J} \Spc_* \,.
    \]
\end{remark}

\subsection{Spin Orientation}

\begin{theorem}
    \label{thm:spinorientation}
    Real and complex K-theory admit orientations by $\E_\infty$-ring homomorphisms
    \begin{equation*}
        M\Spin \to KO \,, \qquad M\Spin^c \to KU \,.
    \end{equation*}
    Moreover, a bundle is $KO$($KU$)-oriented if and only if it admits a $\Spin$($\Spin^c$)-structure.
\end{theorem}
\begin{proof}
    By Proposition~\ref{prop:orientationshomotopies} we must supply an $\E_\infty$-nullhomotopy of 
    \[
        B\Spin \to BO \xrightarrow{J\otimes KO} \Pic(\Mod(KO)) \,.
    \]
    Using the results of Corollary~\ref{cor:spinfiberseq} and Theorem~\ref{thm:spin-orientation-proof},
    this is supplied by commutativity of the following diagram of $\E_\infty$-maps:
     \begin{equation*}
        \begin{tikzcd}
        {B\Spin} && BO && \Pic(\Sp) \\
        {*} && {\Pic(\Bim(\sBan_\R)^{\mathrm{fd}}_\infty)} && {\Pic(\Mod(KO))}
        \arrow[""{name=0, anchor=center, inner sep=0}, from=1-1, to=1-3]
        \arrow[from=1-1, to=2-1]
        \arrow["J", from=1-3, to=1-5]
        \arrow["\Cl"', from=1-3, to=2-3]
        \arrow["{\otimes KO}", from=1-5, to=2-5]
        \arrow[""{name=1, anchor=center, inner sep=0}, from=2-1, to=2-3]
        \arrow["{K^{\gr}}"', from=2-3, to=2-5]
        \arrow[bend right=15,dotted ,from=2-5, to=2-3]
    \end{tikzcd}
    \end{equation*} 
    By the dotted splitting (\cite[Thm.~C]{aretzstehouwer}), existence of a $KO$-orientation is even equivalent to existence of a $\Spin$-orientation.
    The $KU$-version uses the complex Clifford algebra functor $\Cxl$ and the same argument.
\end{proof}

If a vector bundle is not spin, we can still express the $KO$-theory of its Thom space in terms of the graded $K$-theory $K^{\gr}$ of the associated Clifford algebra bundle.
 This provides another Corollary of Theorem~\ref{thm:spin-orientation-proof}.
\begin{corollary}
    Let $V\to X$ be a vector bundle and let $X$ be a finite CW complex.
    Then there is an equivalence of spectra between the $KO$-theory of the Thom space and sections of the $KO$-module bundle $K^{\gr}_{\Cl(-V_x)}$ over $X$:
    \[
        KO(\Th(V))\simeq \lim_{x\in X}K^{\gr}_{\Cl(-V_x)} \simeq \Gamma(X;K^{\gr}_{\Cl(-V)})
    \]
    and there is a natural isomorphism
    \[
        \widetilde{KO}^*(\Th(V)) \cong K^{\gr}_{-*}(\Gamma(X;\Cl(-V)))
    \]
    with the graded $K$-theory $K^{\gr}$ of the Banach algebra $\Gamma(X;\Cl(-V))$.
\end{corollary}
\begin{proof}
    Let $\xi:X\to BO_n$ be a classifying map for $V$.
    Then $\Sigma^\infty \Th(V)=\colim_X J\circ \xi$ and consequently:
    \begin{align*}
        KO(\Th(V))=\map(\colim_X J\circ \xi, KO)& \simeq \lim_X \map_{KO}(J\otimes KO\circ \xi, KO)\\
        &\simeq \lim_X \map_{KO}(K^{\gr}\circ \Cl \circ \xi ,KO) \\
        &\simeq \lim_X K^{\gr}\circ \Cl \circ -\xi \,.
    \end{align*} 
    The fact that taking sections commutes with $K^{\gr}$ follows from~\cite[Prop.~7.1]{hebestreitsagave} or \cite[Prop.~5.30]{bunke2023survey}.
\end{proof}

\subsection{Multiplicativity and slice categories}
\label{sec:multiplicativity}
We will generalize the previous construction to include the tangential structures $p:BH(A;M)\to BO$ of Section~\ref{sec:bimodulebundles}.

Let $V\to X$ be a vector bundle and $\Sigma:\Cl(V-n)\squigrightarrow A$ a bimodule bundle of type $M$.
The pair $(V,\Sigma)$ is classified by a map $\zeta:X\to BH_n(A;M)$ and $V$ is classified by $\xi=p_n\circ \zeta$.
Note that $\Sigma$ defines a natural transformation from $\Cl\circ \xi$ to the constant functor $\const_A$.

Universally, the defining module $M$ defines a homotopy
\begin{equation}
    \label{eq:homotopygeneralH}
    \begin{tikzcd}
    {BH(A;M)} && BO && \Sp \\
    {*} && {\Bim(\sBan_\R)_\infty} && {\Mod(KO)}
    \arrow[""{name=0, anchor=center, inner sep=0}, from=1-1, to=1-3]
    \arrow[from=1-1, to=2-1]
    \arrow["J", from=1-3, to=1-5]
    \arrow["\Cl"', from=1-3, to=2-3]
    \arrow["{\otimes KO}", from=1-5, to=2-5]
    \arrow[""{name=1, anchor=center, inner sep=0}, "A"', from=2-1, to=2-3]
    \arrow["{K^{\gr}}"', from=2-3, to=2-5]
    \arrow[Rightarrow, from=0, to=1, shorten >=.3ex,"M"']
\end{tikzcd}
\end{equation} 
By virtue of Remark~\ref{rmk:mapfromthomspace}, this constructs the natural map \eqref{eq:twistedorientation}.
\begin{theorem}
    \label{thm:multiplicativeorientations}
    Let $A$ be a super Banach algebra and let $M$ be a finitely generated projective $A$-module.
     There are maps of spectra
                \begin{equation}
                    \label{eq:twistedorientation}
                    MH(A;M) \to K^{\gr}_A \,,
                \end{equation}
    satisfying the following:
    \begin{enumerate}[label=(\roman*), ref=(\roman*)]
        \item\label{item:i} Let $MH(A):=\bigoplus_{[M]} MH(A;M)$, where $[M]$ runs over types of finitely generated projective $A$-modules.
                     $MH$ is the left Kan extension of the $J$-homomorphism along the stable Clifford functor $\Cl$.
        \item\label{item:ii} $MH$ defines a lax symmetric monoidal functor $MH:\Bim(\sBan_\R)_\infty\to \Sp$
                     and there is a universal lax symmetric monoidal natural transformation $MH\Rightarrow K^{\gr}$.
    \end{enumerate}
    Therefore, the maps \eqref{eq:twistedorientation} are multiplicative, i.e.\ there are coherently commutative squares
    \begin{equation*}
        \begin{tikzcd}
            MH(A;M)\otimes MH(B;N) \ar[r]\ar[d]& K^{\gr}_{A}\otimes K^{\gr}_B \ar[d]\\
            MH(A\grotimes B;M\grotimes N) \ar[r]& K^{\gr}_{A\grotimes B}
        \end{tikzcd}
    \end{equation*}
    In particular, inserting $A=M=\R$ shows that \eqref{eq:twistedorientation} are $M\Spin\to KO$-module maps.
    \begin{enumerate}[label=(\roman*), ref=(\roman*), start=3]
        \item\label{item:iii} Let $1\le n\le\infty$ and let $A$ be an $\E_n$-algebra in $\Bim(\sBan_\R)_\infty$.
    Then $MH(A)\to K^{\gr}_A$ refines to a map of $\E_n$-rings.
\item\label{item:iv} If $K^{\gr}_A$ is an $\E_1$-ring, e.g.\ as in~\ref{item:iii}, then for every $M$
    the map~\eqref{eq:twistedorientation} is an $H(A;M)$-orientation of $K^{\gr}_A$
    if and only if $[M]\in\pi_0K^{\gr}_A$ is a unit.
    \end{enumerate}
\end{theorem}
\begin{example}
    If $A$ is a graded-commutative super Banach algebra, it defines an $\E_\infty$-algebra in $\Alg(\sBan_\R)_\infty$ and hence in $\Bim(\sBan_\R)_\infty$.
    $MH(A;A)$ is a subring of $MH(A)$ and $MH(A;A)\to K^{\gr}_A$ is an $\E_\infty$-$H(A;A)$-orientation of $K^{\gr}_A$.
\end{example}
\begin{remark}
    For $M=0$, the map \eqref{eq:twistedorientation} is the zero map $MH(A;0)=MO\to K^{\gr}_A$.
\end{remark}
To prove assertion~\ref{item:i}, we will exhibit a universal property of the tangential structures $BH(A;M)$ as a slice category.
\begin{definition}
    Let $\Cl \slice A=BO\times_{\Bim(\sBan_\R)_\infty} (\Bim(\sBan_\R)_\infty)_{/A}$ be the slice category of the functor $\Cl:BO\to \Bim(\sBan_\R)_\infty$ for $A$ a super Banach algebra.
    Explicitly, the objects can be thought of as pairs $(V,\Sigma)$ of virtual rank zero vector spaces $V-\rk(V)$ together with a bimodule $\Sigma:\Cl(V\!-\!\rk V)\squigrightarrow A$.
    A $1$-morphism is a linear isometry $f:V\to W$ and an invertible bimodule map $\phi:\Sigma \to \Sigma'_f$, i.e.\ a $2$-cell $\Sigma\Rightarrow\Sigma'\circ\Cl(f)$.
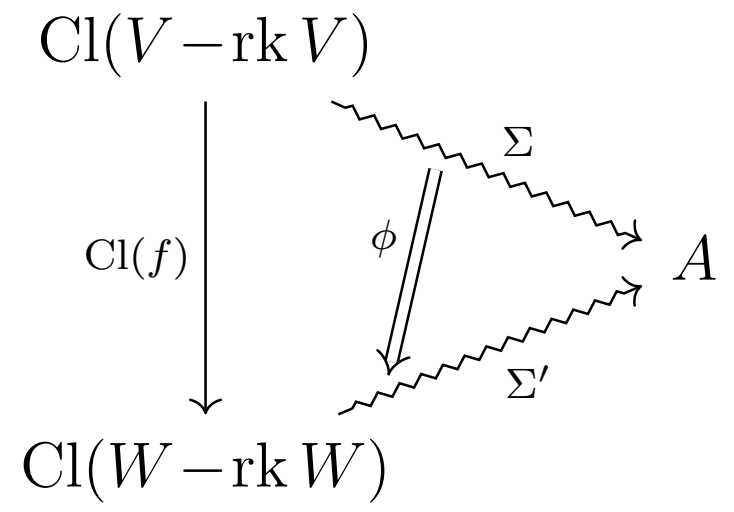
    \begin{equation*}
        \begin{tikzcd}[column sep=large, row sep=normal]
        {\Cl(V\!-\!\rk V)} \\
        & A \\
        {\Cl(W\!-\!\rk W)}
        \arrow["\Sigma", ""{name=0, anchor=center, inner sep=0, pos=0.35}, squiggly, from=1-1, to=2-2]
        \arrow["{\Cl(f)}"', from=1-1, to=3-1]
        \arrow["{\Sigma'}"', ""{name=1, anchor=center, inner sep=0, pos=0.15}, squiggly, from=3-1, to=2-2]
        \arrow["\phi"', Rightarrow, shorten=3pt, from=0, to=1]
    \end{tikzcd}
    \end{equation*}
\end{definition}
The slice category is the universal $\infty$-category with a natural transformation
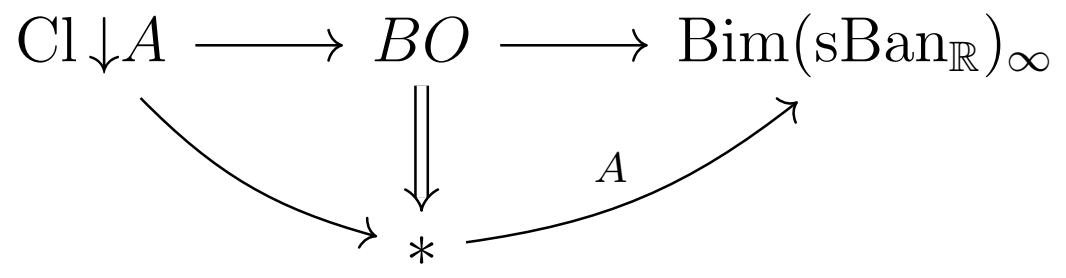
\begin{equation*}
    \begin{tikzcd}
        \Cl\slice A \ar[r] \ar[rd,bend right=15] & BO \ar[r] \ar[d,Rightarrow]& \Bim(\sBan_\R)_\infty \\
        & * \ar[ru,bend right=15,"A"]
    \end{tikzcd}
\end{equation*}
In particular by \eqref{eq:homotopygeneralH}, it receives maps from $BH(A;M)$ for every type $[M]$.
\begin{proposition}
    \label{prop:cliffordslice}
    The slice category $\Cl \slice A$ is an $\infty$-groupoid with homotopy type 
    \begin{equation*}
        \Cl \slice A \simeq \bigsqcup_{M \text{ type}} BH(A;M) \,,
    \end{equation*}
    where $M$ runs over isomorphism classes of finitely generated projective $A$-modules (modulo parity shift).
\end{proposition}
Note that the preceding Proposition works for any super Banach algebra $A$, not necessarily finite-dimensional.
\begin{lemma}
    \label{lem:sliceoverspace}
    Let $X$ be an $\infty$-groupoid and let $F:X\to \calC$ be a functor, $c\in \calC$.
    Then $F\slice c$ is an $\infty$-groupoid and
    \begin{equation*}
        F\slice c \simeq \colim_{x\in X^{\op}} \Map_{\calC}(F(x),c) \,.
    \end{equation*}
\end{lemma}
\begin{proof}
    The slice $\calC_{/c}\to \calC$ is the cartesian fibration classified by $\Map_{\calC}(-,c):\calC^{\op}\to \Spc$.
    The pullback $F\slice c= X\times_{\calC} \calC_{/c}$ is the unstraightening of the composite $\Map_{\calC}(F-,c)$.
    Unstraightening agrees with the colimit.
\end{proof}
    To compute the slice $\Cl\slice A$ we hence have to compute a colimit over $BO$.
    Since $\Cl$ factors over $\tau_{\leq 2}BO$ we can start there.

    The $\E_\infty$-groupoid $\sLine_\R\simeq\Aut_{\Bim(\sBan_\R)_\infty}(\R)$ of real super lines acts on
    $\Map_{\Bim(\sBan_\R)_\infty}(\R,A)\simeq(\Mod^{C_2}_A)^{\cong}$ by precomposition, which is given by the tensor product
    \begin{equation}
        \label{eq:slineaction}
        \begin{aligned}
            \sLine_{\R} \times \Mod^{C_2}_A &\longrightarrow \Mod_A^{C_2} \\
            (L,M) &\longmapsto L\grotimes M \,.
        \end{aligned}
    \end{equation}
    As an $\E_1$-groupoid, $\sLine\simeq \Z/2\times B\{\pm 1\}$ and $B\sLine \simeq \tau_{\leq 2} BO$.
    The following Lemma computes the homotopy quotient of this action.
\begin{lemma}
    \label{lem:slineprincipal}
    For every type $M\neq 0$ in $\Mod^{C_2}_A$, denote the full sub-$\infty$-groupoid spanned by $\{M,\Pi M\}$ by $(\Mod_A^{C_2})^{\cong}_{[M]}$.
    Then
    \begin{equation}
        \label{eq:principal2bundle}
       (\Mod_A^{C_2})^{\cong}_{[M]}\simeq \bigsqcup_{N\in\{M,\Pi M\}/\cong}B\Aut_A(N)
        \longrightarrow B\sAut_A(M)_b
    \end{equation}
    is a principal $\sLine_\R$-``bundle'' for the action~\eqref{eq:slineaction}.
    Here, the left hand side has a single component if $M\cong \Pi M$ and two components otherwise.
    Its classifying map $B\sAut_A(M)_b\to B\sLine_\R\simeq\tau_{\leq2}BO$ is the one of
    the fermionic group $\sAut_A(M)$.
    Equivalently, the map \eqref{eq:principal2bundle} induces an equivalence
    \begin{equation}
        \label{eq:colimitover2type}
        \colim_{\tau_{\leq 2}BO} (\Mod_A^{C_2})^{\cong}_{[M]} \simeq B\sAut_A(M)_b  \,. 
    \end{equation}
    For $M=0$, the action on $(\Mod_A^{C_2})^{\cong}_{[0]}\simeq *$ is trivial and $\colim_{\tau_{\leq 2}BO}(\Mod_A^{C_2})^{\cong}_{[0]}\simeq\tau_{\leq 2}BO$.
\end{lemma}
    The fiber sequence
    \begin{equation*}
        \bigsqcup_{N\in\{M,\Pi M\}/\cong}B\Aut_A(N) \to  B\sAut_A(M)_b \to \tau_{\leq 2}BO  
    \end{equation*}
     is a variant of Proposition~\ref{prop:fiberseqeven}.

\begin{proof}
    We will just prove \eqref{eq:colimitover2type}, as this serves as the definition of a higher principal bundle.
    As an $\E_1$-groupoid $\sLine$ is the product $\Z/2\times B\{\pm 1\}$.
    Model $(\Mod^{C_2}_A)^{\cong}_{[M]}$ by the topological groupoid with
    objects $M,\Pi M$  and even $A$-linear isomorphisms.
    The topological action $2$-groupoid $\mathcal{G}$ associated to the $\sLine$ action has objects $N\in{M,\Pi M}$.
    A 1-morphisms $N\to N'$ is a pair $(a,f)$ where $a\in \Z/2$ and $f:N\to \Pi^a N'$ is a grading-preserving isomorphism of $A$-modules.
    A 2-morphism is $\pm 1 :(a,f)\to (a,-f)$.
    The $2$-groupoid is transitive, i.e.\ there is a single connected component.
    The endomorphisms (or isotropy) of the object $M$ is the topological action groupoid $\{\pm 1\}\ltimes \left(\Aut_A(M)\sqcup \Iso_A(M,\Pi M)\right)$.
    If $M\neq 0$, then $-f\neq f$ for every $f$, so the action is free and the quotient is $\sAut_A(M)/\{\pm 1\}=\sAut_A(M)_b$.
    So the entire topological action $2$-groupoid is equivalent to the topological $1$-groupoid $*//\sAut_A(M)_b$. 
    We obtain
    \begin{equation*}
        \colim_{\tau_{\leq 2}BO} (\Mod^{C_2}_A)^{\cong}_{[M]}  = |\mathcal{G}| \simeq B\sAut_A(M)_b  \,.
    \end{equation*}
    It is easy to see that these are equivalences over $\tau_{\leq 2}BO$.
    If $M=0$, then $(\Mod_A^{C_2})_{[M]}^{\cong} = *$ and $\colim_{\tau_{\leq 2}BO }*=\tau_{\leq 2}BO$.
    This yields the claim.
\end{proof}
\begin{proof}[Proof of Proposition~\ref{prop:cliffordslice}]
    By Lemma~\ref{lem:sliceoverspace}, the slice can be expressed as a colimit over $BO$, i.e.\ as a homotopy quotient:
    \begin{equation*}
        \Cl\slice A \simeq \colim_{BO}\Map_{\Bim(\sBan_\R)_\infty}(\Cl(-),A) \,.
    \end{equation*}
    Since the stable Clifford functor $\Cl$ takes value at the tensor unit, the mapping space is just 
    \[
        \Map_{\Bim(\sBan_\R)_\infty}(\R,A)\simeq (\Mod_A^{C_2})^{\cong} \simeq \bigsqcup_{M \text{ fgp}} B\Aut_A(M)
    \]
    By Proposition~\ref{prop:cliffordpostnikov}, the precomposition action of $BO$ factors over $\tau_{\leq 2}BO=B\sLine_\R$.
    The action of $\sLine$ on $\Mod^{C_2}_A$ is then given by tensor product and using Lemma~\ref{lem:slineprincipal}, the homotopy quotient is
    
    \begin{equation*}
        \colim_{\tau_{\leq 2} BO}  \bigsqcup_{M \text{ fgp}} B\Aut_A(M) = \tau_{\leq 2}BO\sqcup\bigsqcup_{M\neq 0 \text{ type}} B\sAut_A(M)_b\,,
    \end{equation*}
    where the first summand comes from $M=0$.
    Since colimits of spaces are universal (stable under pullback), we have
    \begin{equation*}
        \Cl\slice A \simeq BO\sqcup\bigsqcup_{M\neq 0 \text{ type}} BO\times_{\tau_{\leq 2} BO} B\sAut_A(M)_b \,.
    \end{equation*}
    For $M\neq 0$ we are done by the definition of the fermionic tensor product $B\sAut_A(M)_b\times_{\tau_{\leq 2} BO} BO \simeq B(\Pin^+ \grotimes \sAut_A(M)^{\op})_{\ev}$, which reproduces $BH(A;M)$.
    For $M=0$, the summand is $BO=BH(A;0)$, since $\Aut(\R^n,\Cl_n\grotimes 0)=O_n$; its Thom spectrum is $MH(A;0)=MO$.
\end{proof}

We are now ready to prove multiplicativity.
\begin{proof}[Proof of Theorem~\ref{thm:multiplicativeorientations}]
    The stable Clifford functor $\Cl\colon BO\to\Bim(\sBan_\R)_\infty$ is symmetric monoidal.
    By Lemma~\ref{lem:kanextensionmonoidal}, applied to $\gamma=\Cl$ and $\calE=\Mod(KO)$, the lax symmetric monoidal natural equivalence
    $\alpha\colon J\otimes KO\simeq K^{\gr}\circ\Cl$ of Theorem~\ref{thm:spin-orientation-proof} is adjoint to a lax symmetric monoidal transformation
    \[
        \Lan_{\Cl}(J\otimes KO)\Longrightarrow K^{\gr}\,.
    \]
    Since $-\otimes KO$ preserves colimits, $\Lan_{\Cl}(J\otimes KO)\simeq(\Lan_{\Cl} J)\otimes KO$.
    The pointwise formula for left Kan extensions \cite[Sec.~4.3.3]{HTT} computes $\Lan_{\Cl} J$ at $A$ as a colimit over the slice category $BO\times_{\Bim(\sBan_\R)_\infty}(\Bim(\sBan_\R)_\infty)_{/A}=\Cl\slice A$.
    Writing $\pi_A\colon\Cl\slice A\to BO$ for the projection, Proposition~\ref{prop:cliffordslice} gives
    \[
        (\Lan_{\Cl} J)(A)\simeq\colim_{\Cl\slice A}J\circ\pi_A\simeq\colim_{\bigsqcup_{[M]}BH(A;M)}J\simeq\bigoplus_{[M]}MH(A;M)=MH(A)\,.
    \]
    This proves \ref{item:i}.
    The composite lax symmetric monoidal natural transformation
    \[
        MH\Rightarrow MH\otimes KO \simeq \Lan_{\Cl}(J\otimes KO) \Rightarrow K^{\gr}
    \] 
    proves \ref{item:ii}.
    It is easily seen that the map out of the component $MH(A;M)$ is induced from the cocone \eqref{eq:homotopygeneralH}.

    To prove \ref{item:iv}, we consider the natural transformation
    \begin{equation*}
        \begin{tikzcd}
            BH(A;M)\ar[r]\ar[rd,bend right=15]& BO \ar[r,"J\otimes K^{\gr}_A"]\ar[d,Rightarrow] & \Mod(K^{\gr}_A)\\
            & * \ar[ru,bend right=15,"K^{\gr}_A"'] &
        \end{tikzcd}
    \end{equation*}
    obtained from \eqref{eq:homotopygeneralH} by $K^{\gr}_A$-linearizing.
    The map $MH(A;M)\to K^{\gr}_A$ is an orientation if and only if this natural transformation is an equivalence.
    Equivalence can be checked at the unique point $*\to BH(A;M)$, where the natural transformation is multiplication by $[M]\in \pi_0K^{\gr}_A$.

    For~\ref{item:iii}, lax symmetric monoidal functors and transformations preserve $\E_n$-algebras and maps between them.
    Hence $MH(A)\to K^{\gr}_A$ is a map of $\E_n$-rings.
\end{proof}

\begin{lemma}
    \label{lem:kanextensionmonoidal}
    Let $\gamma:\calC \to \calD$ be a symmetric monoidal functor between small symmetric monoidal $\infty$-categories. Let $\calE$ be presentably symmetric monoidal.
    Then left Kan extension refines to an adjunction
    \begin{equation*}
        \begin{tikzcd}
            \gamma^*:\Fun^{\mathrm{lax}}(\calD,\calE) \ar[r, shift right=.5ex] & \Fun^{\mathrm{lax}}(\calC,\calE) \ar[l, shift right=.5ex]: \Lan_\gamma
        \end{tikzcd}
    \end{equation*}
    on the $\infty$-categories of lax symmetric monoidal functors and lax symmetric monoidal natural transformations.
\end{lemma}
\begin{proof}
    By \cite[Cor.~3.8]{Nikolaus2016stable}, left Kan extension is a symmetric monoidal left adjoint functor for the Day convolution symmetric monoidal structure on functor $\infty$-categories and $\gamma^*$ is lax symmetric monoidal.
    By \cite[Rem.~7.3.2.13]{HA} the adjunction restricts to the subcategories of commutative algebra objects, which are lax symmetric monoidal functors.
\end{proof}

\begin{remark}[Size issues]
    \label{rmk:size}
    In the proof of Theorem~\ref{thm:multiplicativeorientations}, Lemma~\ref{lem:kanextensionmonoidal} is applied to $\Cl\colon BO\to\Bim(\sBan_\R)_\infty$, whose target is not small.
    Passing to a larger universe as in \cite[Rem.~2.8]{Nikolaus2016stable}, the lemma applies with target the $\infty$-category of large $KO$-modules, which contains $\Mod(KO)$ as a full symmetric monoidal subcategory closed under small colimits.
    The Kan extension lands in $\Mod(KO)$ and agrees with the one computed there, since its value at $A$ is a colimit over $\Cl\slice A$, which is essentially small by Lemma~\ref{lem:sliceoverspace} and~\eqref{eq:mappingspacesbimban1}.
    Alternatively, one can restrict to an arbitrary small subcategory of the target.
\end{remark}

\subsection{Examples}
\label{sec:examples}
\subsubsection{The tenfold way}
In Section~\ref{sec:equivalencesandirreps}, we showed that we can recover the spacetime structure groups $H_n(s)$ of \cite{FreedHopkins} from super division algebras as the structure groups $H_n(D;D)$ classifying pairs $(V,\Sigma)$ of a vector bundle $V$ and a bimodule bundle $\Sigma:\Cl(V\!-\!\rk V)\squigrightarrow D$ with irreducible fibers.
By Theorem~\ref{thm:multiplicativeorientations}, there are $M\Spin$-module maps
\begin{equation*}
    MH(D;D) \to K^{\gr}_D
\end{equation*}
Using that $K^{\gr}_{D}$ is a shift of $KO$ (or $KU$ for $D$ a complex super division algebra), we finally obtain $M\Spin$-module maps
\begin{equation}
    \label{eq:twistedABSFH}
    MH(s) \to \Sigma^s KO \,, \qquad MH^c(s) \to \Sigma^{s} KU \,,
\end{equation}
which reproduces the maps of \cite[(9.42)]{FreedHopkins}.

\begin{lemma}[{\cite[(10.2)--(10.3)]{FreedHopkins}}]
    The Thom spectra split as $M\Spin$-modules:
    \[
        MH(s)\simeq M\Spin\wedge X(s),\qquad
        MH^c(s)\simeq M\Spin\wedge X^c(s),
    \]
    with factors as in~\eqref{eq:100} and~\eqref{eq:101}.
    For the tangential variants $MTX=\colim_X(-J)$ of a stable structure map $X\to BO$, we have $MTH(s)\simeq MH(-s)$.
\end{lemma}
By $M\Spin$-linearity the maps \eqref{eq:twistedABSFH} factor as
\[
   MH(s) \simeq M\Spin\wedge X(s)\longrightarrow ko \wedge X(s)
    \longrightarrow\Sigma^{s}KO.
\]
\begin{proof}
    Under the stable splitting of Lemma~\ref{lem:pullbackfermionicgrps} $BH(s)\simeq B\Spin \times BO(D(s))_b$ compatible with $J$.
    Therefore
     \[
       MH(s)= \colim_{BH(s)} J \simeq \colim_{B\Spin} J \wedge \colim_{BO(D(s))_b} J\simeq M\Spin \wedge X(s) \,. \qedhere
    \]
\end{proof}

We use the group labels and Clifford types of
\cite[Tables~(9.24)--(9.25)]{FreedHopkins}.
See also \cite[Fig.~31]{beaudrycampbell}.
\begin{equation}
    \label{eq:100}
    \begingroup
    \setlength{\tabcolsep}{5pt}
    \begin{tabular}{c l c c c c c}
    \toprule
    $s$&$H$&$K$&Cartan&$D$&$O(D)_b$&$X(s)$\\ \midrule
    $0$&$\Spin$&$O_1$&D&$\R$&$*$&$\mathbb S$\\
    $-1$&$\Pin^+$&$O_1$&DIII&$\Cl_{-1}$&$O_1$&$\Sigma^{-1}MO_1$\\
    $-2$&$\Pin^{\tilde c+}$&$U_1$&AII&$\Cl_{-2}$&$O_2$&$\Sigma^{-2}MO_2$\\
    $-3$&$\Pin^{h-}$&$Sp_1$&CII&$\Cl_{-3}$&$O_3$&$\Sigma^{-3}MO_3$\\
    $4$&$\Spin^h$&$Sp_1$&C&$\mathbb H$&$SO_3$&$\Sigma^{-3}MSO_3$\\
    $3$&$\Pin^{h+}$&$Sp_1$&CI&$\Cl_{+3}$&$O_3$&$\Sigma^3 MTO_3$\\
    $2$&$\Pin^{\tilde c-}$&$U_1$&AI&$\Cl_{+2}$&$O_2$&$\Sigma^2 MTO_2$\\
    $1$&$\Pin^-$&$O_1$&BDI&$\Cl_{+1}$&$O_1$&$\Sigma MTO_1$\\
    \bottomrule
    \end{tabular}
    \endgroup
\end{equation}
\begin{equation}
    \label{eq:101}
    \begingroup
    \setlength{\tabcolsep}{5pt}
    \begin{tabular}{c l c c c c c}
    \toprule
    $s$&$H^c$&$K$&Cartan&$D^c$&$O(D^c)_b$&$X^c(s)$\\ \midrule
    $0$&$\Spin^c$&$U_1$&A&$\C$&$U_1$&$\Sigma^{-2}MU_1$\\
    $1$&$\Pin^c$&$U_1$&AIII&$\Cxl_{-1}$&$U_1\times O_1$
    &$\Sigma^{-3}(MU_1\wedge MO_1)$\\
    \bottomrule
    \end{tabular}
    \endgroup
\end{equation}

We will now make some remarks on potential physical relevance.
\begin{remark}
    Physically, free fermion phases in dimension $n$ with $H(s)$-symmetry are classified by $
    KO$-theory and unitary interacting phases with $H(s)$-symmetry are postulated to be classified by the Anderson dual $I_{\Z}MTH(s)$.

    Roughly speaking, a free fermion particle corresponds to an irreducible representation on a super Hilbert space.
    By Schur's Lemma, there are ten classes of these labelled by the ten super division algebras.
    Free fermion phases in spacetime dimension $n$ with symmetry type $s$ are classified by $\pi_{3-n-s}KO$~\cite[Thm.~9.63]{FreedHopkins} and they decompose into symmetry types.
    Mathematically, this corresponds to the fact that representations are Morita invariant. 
    Stehouwer~\cite{stehouwermorita} argues that symmetries of interacting phases do not have such a tenfold classification and that the spacetime group of a symmetry type is ambiguous.
    In our language, he analyzes the example $s=0$ with algebra $A=\R$ and the two different modules $M=\R$ or $M=\R^2$.
\end{remark}
        
\begin{remark}
    In the previous remark we argued that to a symmetry type $s$ there belong different symmetry groups and the periodic assignment $s\mapsto H(s)$ as in \cite{FreedHopkins} is not unique.
    We now outline a multiplicativity deficit of the periodicity in $s$, namely only a partial multiplicativity.
    The multiplicative structure is
    \begin{equation*}
        MH(\Cl_s;S) \otimes MH(\Cl_t;T) \to MH(\Cl_s\grotimes\Cl_t; S\grotimes T) \,.
    \end{equation*}
    Suppose $S,T$ are irreducible.
    For $\mod 8$-indices $s,t\in\{-3,\ldots,4\}$, a dimension argument shows that $S\grotimes T$ is irreducible precisely in the following cases, up to exchanging $s$ and $t$:
    \begin{equation*}
        \begin{aligned}
            (s,t)&=(0,u)\quad\text{for }u\in\{-3,\ldots,4\},\\
            \text{or}\quad(s,t)&\in\{(-2,-1),(-1,-1),(-1,4),(1,1),(1,2),(1,4)\}.
        \end{aligned}
    \end{equation*}
    This gives a partial multiplicativity, e.g.\ for $(s,t)=(1,4)$
    \begin{equation*}
        M\Pin^- \otimes M\Spin^h \to M\Pin^{h-} \,.
    \end{equation*}
    If $S,T$ are the left regular module, this is always multiplicative, but then $MH(\Cl_s;S)\simeq MH(s)$ only for $-3\leq s\leq 3$.    
\end{remark}

\begin{remark}
    The (unstable) multiplicativity $MTH_n(A;M)\otimes MTH_m(B;N) \to MTH_{n+m}(A\grotimes B;M\grotimes N)$ is related to Kaluza--Klein compactification of TQFTs as discussed in \cite[Sec.~6 \& Sec.~7.3]{yamashitayonekura}.
\end{remark}
\subsubsection{Unstable variants}
    \label{rmk:ebert2}    
    We now show that by the methods of Section~\ref{sec:multiplicativity} \eqref{eq:homotopygeneralH} it is also possible to recover certain unstable maps from the literature.
    This is a continuation of Remark~\ref{rmk:ebert1}
    
    For $A$ a Real super $C^*$-algebra, Ebert's topological index~\cite[Definition~2.22]{EbertIndex} is a map of spectra
    \begin{equation}
        \label{eq:ebertmap}
        MT\theta_A(d)\to K^{\gr}_{A}\,,
    \end{equation}
    where $MT\theta_A(d)=\bigoplus_{N} MTH_d(A^{\op}\grotimes \Cl_{-d};N)$.
    The finitely generated projective $A^{\op}\grotimes \Cl_{-d}$-module $N$ has a dual $A\grotimes \Cl_d$-module $N^\vee$, which may be viewed as a morphism $N^\vee:\Cl_{-d}\squigrightarrow A$.

    \[
    \begin{tikzcd}
    {BH_d(A^{\op} \grotimes \Cl_{-d};N)} && BO_d && \Sp \\
    {*} && {\Bim(\sBan_\R)_\infty} && {\Mod(KO)}
    \arrow[""{name=0, anchor=center, inner sep=0, pos=0.4}, from=1-1, to=1-3]
    \arrow[from=1-1, to=2-1]
    \arrow["-J", from=1-3, to=1-5]
    \arrow["\Cl^{\op}"', from=1-3, to=2-3]
    \arrow["-\otimes KO", from=1-5, to=2-5]
    \arrow[""{name=1, anchor=center, inner sep=0, pos=0.6}, "A"', from=2-1, to=2-3]
    \arrow["K^{\gr}"', from=2-3, to=2-5]
    \arrow[Rightarrow, from=0, to=1, shorten >=.3ex,"N^\vee"']
    \end{tikzcd}
    \]
    The above diagram replaces \eqref{eq:homotopygeneralH}.
    Comparing $K$-theory conventions of \cite{EbertIndex} and \cite{aretzstehouwer}, it is easy to see that the Thom classes are the same and that this indeed reproduces Ebert's topological index~\eqref{eq:ebertmap}.

\subsubsection{Group algebras and assembly}
\label{sec:groupalgebras}
Let $G$ be a finite fermionic group with central element $c$.
The fermionic group algebra is the finite dimensional semisimple super algebra $\R^f[G]=\R[G]/(c=-1)$.
$G\mapsto \R^f[G]$ is the symmetric monoidal left adjoint to the fermionic unit group functor $A\mapsto A^\times$.
The adjunction unit is the inclusion $G\to \R^f[G]^\times$.
Let $H_n^G=(\Pin_n^+\grotimes G)_{ev}$.
The adjunction unit induces $H_n^G\to H(\R^f[G];\R^f[G])$ and Theorem~\ref{thm:multiplicativeorientations} supplies maps
\begin{equation}
    \label{eq:grouporientation}
    MH^G\to MH(\R^f[G];\R^f[G])\to K^{\gr}_{\R^f[G]}
\end{equation}
These maps are multiplicative in $G$.
In particular, if $G$ is graded commutative, then \eqref{eq:grouporientation} is an $\E_\infty$-orientation of $K^{\gr}_{\R^f[G]}$.
The right hand side is the $(w_1(G),w_2(G))$-twisted $G_b$-equivariant $KO$-theory of a point.

For purely even infinite discrete groups $\Gamma$, we can instead consider $\ell^1(\Gamma)$ and the trivial extension to a fermionic group $\hat{\Gamma}=\Gamma\times \Z/2$.
In this case, there is a splitting $H^{\hat{\Gamma}}_n= \Spin_n\times \Gamma$ and hence we obtain a map
\[
    M\Spin \wedge B\Gamma_+ \to KO\wedge B\Gamma_+ \to  KO_{\ell^1(\Gamma)} \,.
\]
The second map is a real Banach (Bost) assembly map.
There is a further map $KO_{\ell^1(\Gamma)} \to KO_{C^*_{\R}\Gamma}$ for the maximal and reduced group $C^*$-algebra.
An index theorem links the composite to the index of the Dirac operator twisted by the Mishchenko--Fomenko bundle~\cite[Rmk.~7.3]{LNS} \cite[Ex.~2.13]{EbertIndex}.

\subsubsection{Sullivan orientation}
\label{sec:sullivan}
This section sketches how to obtain $\E_\infty$-ring homomorphisms $MSO\to KO[\tfrac{1}{2}]$ from the Canonical-Anticommutation-Relations algebra.
Consider the inclusions
\begin{equation*}
    M_{2^n}(\R)\rightarrow M_{2^{n+1}}(\R)\,, \qquad T\mapsto T\oplus T \,,
\end{equation*}
where the ungraded Banach algebras carry the unique $C^*$-norm and let $A=\colim_{\N} M_{2^n}(\R)$, the Banach completion of the union in the supremum norm.
This is a special case of a UHF algebra, the CAR algebra.
By \cite[Thm.~B]{aretzstehouwer}, $K^{\gr}$ preserves filtered colimits along contractive homomorphisms. The transition maps are isometric and act as multiplication by two, so
\begin{equation*}
    K_{A}=\colim_{\N} \left(KO\xrightarrow{\cdot 2}KO \rightarrow\dots \right) =KO[\tfrac{1}{2}] \,.
\end{equation*}
Note that we may write $M_{2^n}(\R)=M_{2}(\R)^{\otimes n}$ and $A=\bigotimes_{i=1}^{\infty}M_2(\R)$ with the transition maps becoming $M_2(\R)\otimes \R \to M_2(\R)\otimes M_2(\R)$.
The algebra $A$ carries an $\E_\infty$-multiplication via
\begin{align*}
    M_2(\R)^{\otimes \infty} \otimes M_2(\R)^{\otimes \infty} &\to M_2(\R)^{\otimes \infty} \\
    (a_1\otimes a_2\otimes \dots,b_1\otimes b_2\otimes \dots) &\mapsto (a_1\otimes b_1\otimes a_2\otimes b_2\otimes\dots ) \,.
\end{align*}
\begin{remark}
    \cite[Ex.~1.14]{TomsWinter} proves that $A$ is a strongly self-absorbing $C^*$-algebra, i.e.\ that the unit $\eta:1\to A$ induces an equivalence $\id \otimes \eta:A\to A\otimes_{\min} A$ in the $\infty$-category of $C^*$-algebras.
    This is called an idempotent object in~\cite[Def.~4.8.2.1]{HA} and idempotent objects carry a natural $\E_\infty$-multiplication inverse to $\id \otimes \eta$~\cite[Prop~4.8.2.9]{HA}.
    The natural comparison map $A\otimes_\pi A\to A\otimes_{\min}A$ supplies $A$ with an $\E_\infty$-structure in $(\Bim(\sBan_{\R}),\otimes_\pi)$.
\end{remark}
By Theorem~\ref{thm:multiplicativeorientations}, there is an $\E_\infty$-orientation
\begin{equation}
    \label{eq:universal12orientation}
    MH(A;A)\to KO[\tfrac{1}{2}] \,.
\end{equation}
We have $H_r=H_r(A;A)= (\Pin_r^+\grotimes A^\times)_{\ev}= \Spin_r\times A^\times/\langle(-1,-1) \rangle$ fitting into a short exact sequence
\begin{equation*}
    1\to A^\times \to H_r \to SO(r) \to 1 \,.
\end{equation*}
There is an explicit splitting using Clifford algebras:
Let $\tilde{g}\in \Spin(r)$ be a lift of $g\in SO(r)$.
The standard left action $\rho_r$ of $\Spin(r)$ on $\lvert\Cl_r\rvert\cong \R^{2^r}$ supplies a splitting
\begin{equation}
    \label{eq:splittingmap}
    \begin{aligned}
         SO(r) &\rightarrow H_r \\
    g & \mapsto [\tilde{g},\rho_r(\tilde{g})] \,.
    \end{aligned}
\end{equation}
This is compatible with stabilization.
We conjecture that this lifts to an $\E_\infty$-splitting.
In particular, there is also a (non-canonical) splitting $BH(A;A)\simeq BSO\times BA^\times$.
A splitting can also be obtained by obstruction theory.
Furthermore, splittings of the fiber sequence $BA^\times\to BH(A;A)\to BSO$ are in bijection with orientations
\begin{equation}
    \label{eq:msoorientation}
    MSO\to KO[\tfrac{1}{2}] \,,
\end{equation}
which map $1\in \pi_0 MSO\cong \Z$ to $1\in \pi_0KO[\tfrac{1}{2}]\cong \Z[\tfrac{1}{2}]$.
For multiplicative orientations this condition is automatic. We state the $\E_\infty$-version:
\begin{proposition}
    The space of $\E_\infty$-sections of $p:BH(A;A)\to BSO$ is equivalent to the space of $\E_\infty$-orientations $MSO\to KO[\tfrac{1}{2}].$
    Equivalently, $BH(A;A)\simeq \fib(BSO\to BSL_1(KO[\tfrac{1}{2}]))$.
\end{proposition}
\begin{proof}[Proof sketch.]
    Write $j:BSO\to BSL_1(KO[\tfrac{1}{2}])=\tau_{\geq 2} BGL_1(KO[\tfrac{1}{2}])$ for the $J$-homomorphism.
    The space of $\E_\infty$-orientations is the space of $\E_\infty$-sections of $\fib(j)\to BSO$.
    The orientation~\eqref{eq:universal12orientation} provides a map $BH(A;A)\to \fib(j)$ over $BSO$.
    To show that this is an equivalence we consider the induced map on fiber sequences: 
    \[
    \begin{tikzcd}
        BA^{\times} \ar[r]\ar[d]&BH(A;A) \ar[d]\ar[r]& BSO \ar[d,equals]\\
        SL_1(KO[\tfrac{1}{2}]) \ar[r]& \fib(j) \ar[r]& BSO
    \end{tikzcd}
    \]
    We show that the map is an equivalence on fibers and hence on total spaces.
    Now $\Omega^\infty KO[\tfrac{1}{2}]=\Omega^\infty K_A=K_0(A)\times B\GL_\infty(A)$.
    We also have that $A^\times \to \GL_\infty(A)$ is an equivalence.
    Hence the underlying space of $\GL_1(KO[\tfrac{1}{2}])$ is $\Z[\tfrac{1}{2}]^\times \times BA^\times$ and $BA^\times \to \mathrm{SL}_1(KO[\tfrac{1}{2}])$ is an equivalence.
\end{proof}
\begin{remark}    
    There are many different $SO$-orientations of $KO[\tfrac{1}{2}]$.
    An $\E_\infty$-refinement of the signature $MSO\to KO[\tfrac{1}{2}]$ is constructed in \cite[Thm.~E]{LNS}; see also \cite{banagl} for a detailed discussion and the relation to the work of Sullivan and Ranicki.
    Its characteristic series is $\frac{x/2}{\tanh(x/2)}$.
    By contrast, the composite $MSO\to MSO[\tfrac{1}{2}]\simeq M\Spin[\tfrac{1}{2}]\to KO[\tfrac{1}{2}]$ induced by the $\Spin$-orientation has characteristic series $\frac{x/2}{\sinh(x/2)}$ (Appendix~\ref{app:comparison}).
    The methods of Appendix~\ref{app:comparison} also show that the explicit splitting~\eqref{eq:splittingmap} has characteristic series $\frac{x}{\sinh x}$.
\end{remark}

\newpage
\appendix
\section{Lax monoidal mapping spaces and duality}
\label{app:mappingspaces}

\begin{lemma}
    \label{lem:mappingspaces}
    Let $\calC$ and $\calD$ be symmetric monoidal $\infty$-categories.
    \begin{enumerate}[label=(\arabic*), ref=(\arabic*)]
        \item\label{item:Maplax} The functor $\Map_{\calC}: \calC^{\op}\times \calC \to \Spc$ is lax symmetric monoidal.
        \item\label{item:compatibilityop}  There is a lax symmetric monoidal equivalence
  $\Map_{\calC^{\op}}(c,c')\simeq\Map_{\calC}(c',c)$ of functors $\calC\times\calC^{\op}\to\Spc$.
        \item\label{item:Mapnatural} If $F:\calC\to \calD$ is symmetric monoidal, then $F:\Map_\calC\to \Map_{\calD}\circ (F^{\op}\times F)$ is a lax symmetric monoidal natural transformation.
        \item\label{item:dualityfunctor} If every object of $\calD$ is dualizable, there is a symmetric monoidal functor $(-)^\vee:\calD^{\op}\to \calD$ which assigns to every object a dual and a symmetric monoidal natural equivalence $d\to (d^{\vee})^{\vee}$.
        \item\label{item:dualityadjunction} If every object of $\calD$ is dualizable, there is a lax symmetric monoidal natural equivalence of functors $\calD^{\op}\times\calD\to\Spc$
        \begin{equation*}
            \Map_{\calD}(d,d')\xrightarrow{\simeq} \Map_{\calD}(\mathbbm{1},d^\vee\otimes d') \,.
        \end{equation*}
    \end{enumerate}
\end{lemma}
A version of \ref{item:dualityfunctor} is proven in \cite[Thm.~5.21]{heinelopezspitzweck}, but it does not suffice to prove \ref{item:dualityadjunction}.
\begin{proof}
    Lurie defines the twisted arrow category $\TwArr(\calC)$ by the limit-preserving functor on $\sSet$, mapping $\calC$ to the simplicial set $[n] \mapsto \calC([n]\ast [n]^{\op})$~\cite[5.2.1.1]{HA}, where for $I,J$ linearly ordered $I\ast J$ is the linear order on $I\sqcup J$ in which $I<J$.
    The natural maps $I\to I\ast I^{\op} \leftarrow I^{\op}$ provide maps $\calC^{\op} \leftarrow \TwArr(\calC) \to \calC$.
    The twisted arrow category $\TwArr(\calC)\to \calC^{\op} \times \calC$ is a right fibration and it is classified by the functor $\Map_{\calC}:\calC^{\op}\times \calC \to \Spc$~\cite[5.2.1.11]{HA}.
    $\TwArr$ restricts to an endofunctor on $\infty$-categories and preserves products. It hence preserves symmetric monoidal $\infty$-categories. 
    Furthermore, the canonical isomorphisms $I\cong I^{\op}$ induce a natural isomorphism $\TwArr(\calC)\cong \TwArr(\calC^{\op})$ of simplicial sets or $\infty$-categories.
    The monoidal version of straightening~\cite[Thm.~B]{ramzimonoidal} yields a lax symmetric monoidal structure on $\Map_{\calC}:\calC^{\op}\times \calC\to \Spc$ and a lax symmetric monoidal equivalence $\Map_{\calC}(c,c')\simeq \Map_{\calC^{\op}}(c',c)$.
    This proves \ref{item:Maplax} and \ref{item:compatibilityop}.
    A symmetric monoidal functor $F:\calC \to \calD$ induces a symmetric monoidal $F_*:\TwArr(\calC)\to \TwArr(\calD)$ and by straightening a monoidal natural transformation as in \ref{item:Mapnatural}. 

    Enlarging the universe, we may assume that $\calC$ and $\calD$ are small.
    Equip $\mathcal{P}(\calC):=\Fun(\calC^{\op},\Spc)$ with the Day convolution product $\circledast$~\cite[Ex.~2.2.6.17]{HA}.
    By its universal property~\cite[Def.~2.2.6.1, Ex.~2.2.6.9]{HA}, for every symmetric monoidal $\infty$-category $\calE$ there is an equivalence
    \begin{equation}
        \label{eq:dayadjunction}
        \Fun^{\mathrm{lax}}(\calE\times \calC^{\op},\Spc)\simeq \Fun^{\mathrm{lax}}(\calE,\mathcal{P}(\calC))
    \end{equation}
    which on underlying functors is currying.
    For $\calE=\calC$ the lax symmetric monoidal structure on $\Map$ corresponds to a lax symmetric monoidal structure on the Yoneda embedding $j_\calC:\calC\to \mathcal{P}(\calC)$.      
    The Yoneda embedding is furthermore strong symmetric monoidal: By \cite[Ex.~2.2.6.17]{HA},
    Day convolution is left Kan extension along $\otimes$ and hence $j_{\calC}(x)\circledast j_{\calC}(y) = \otimes_!(j_{\calC}(x)\boxtimes j_{\calC}(y))\simeq \otimes_!(j_{\calC\times \calC}(x,y))=j_{\calC}(x\otimes y)$ and the structure maps here match with those obtained from currying.

    The functor 
    \begin{align}
        \label{eq:fiberfunctor}
        \calD \times \calD \xrightarrow{\otimes} \calD \xrightarrow{\Map(\mathbbm{1},-)} \Spc
    \end{align}
    is lax symmetric monoidal and adjoint under \eqref{eq:dayadjunction} to a lax symmetric monoidal functor $\calD\to \Fun(\calD,\Spc), d\mapsto \Map(\mathbbm{1},d\otimes -)$.
    This functor is symmetric monoidal:
    We claim that the monoidality map
    \begin{equation}
        \label{eq:monoidalitymap}
        \Map(\mathbbm{1},d\otimes - )\circledast \Map(\mathbbm{1},d'\otimes -)\to \Map(\mathbbm{1},(d\otimes d')\otimes -)
    \end{equation}
    is an equivalence.
    The monoidality map \eqref{eq:monoidalitymap} is induced from the tensor product composed with the symmetric braiding
    \begin{equation*}
             \Map(\mathbbm{1},d\otimes c)\times \Map(\mathbbm{1},d'\otimes c')\xrightarrow{\otimes}  \Map(\mathbbm{1},d\otimes c\otimes d'\otimes c')\simeq  \Map(\mathbbm{1},d\otimes d'\otimes c\otimes c')\,,
    \end{equation*}
    naturally in $c,c'\in \calD$.
    Applying this to duals $c=d^\vee$ and $c'=d'^\vee$ and their coevaluations $\mathbbm{1}\to d\otimes d^\vee$, the coevaluation $\coev_d\otimes \coev_{d'}$ is in the image of \eqref{eq:monoidalitymap} and it is a coevaluation for $d\otimes d'$.
    We will show below that this implies that \eqref{eq:monoidalitymap} is an equivalence.

    Again using that $\calD$ has duals, the functor~\ref{eq:fiberfunctor} factors through the symmetric monoidal full subcategory of corepresentable functors, which is equivalent to $\calD^{\op}$.
    This exhibits a lax symmetric monoidal dual functor $(-)^\vee:\calD\to \calD^{\op}$.
    Under \eqref{eq:dayadjunction} this is adjoint to a monoidal equivalence $\Map_{\calD}(\mathbbm{1},d \otimes d')\simeq \Map_{\calD}(d^{\vee},d')$ of
    lax symmetric monoidal functors $\calD\times \calD\to \Spc$.
    In the following diagram the identity is in the image of the top horizontal composite since the image of \eqref{eq:monoidalitymap} contains a coevaluation.
    \begin{equation*}
        \begin{tikzcd}[column sep=1.6em]
            \Map(d^\vee\otimes d'^\vee,-) \arrow[r, phantom, "\simeq"]
            & \Map(d^\vee,-)\circledast \Map(d'^\vee,-) \arrow[r] \arrow[d, "\simeq"']
            & \Map((d\otimes d')^\vee,-) \arrow[d, "\simeq"] \\
            & \Map(\mathbbm{1},d\otimes -)\circledast \Map(\mathbbm{1},d'\otimes -) \arrow[r, "\eqref{eq:monoidalitymap}"]
            & \Map(\mathbbm{1},d\otimes d'\otimes -)
        \end{tikzcd}
    \end{equation*}
    Therefore, the top composite is an equivalence and hence also \eqref{eq:monoidalitymap}. 

    The chain of monoidal equivalences of lax symmetric monoidal functors $\calD \times \calD \to \Spc$
    \begin{align*}
         \Map_{\calD}(d^\vee,d') &\simeq \Map_{\calD}(\mathbbm{1},d\otimes d')\simeq \Map_{\calD}(\mathbbm{1},d'\otimes d)\\
         &\simeq \Map_{\calD}(d'^{\vee},d)\simeq \Map_{\calD^{\op}}(d,(d')^{\vee\op})
    \end{align*}
    induces a monoidal adjunction between $(-)^\vee:\calD\to \calD^{\op}$ and $(-)^{\vee\op}:\calD^{\op}\to \calD$,
     which is moreover an equivalence: its unit and counit are the canonical map to the double dual, which is an equivalence.
    This proves \ref{item:dualityfunctor}.
    We obtain a monoidal equivalence of lax symmetric monoidal functors $\calD^{\op}\times \calD\to \Spc$:
    \begin{equation*}
        \Map_{\calD}(d,d')\simeq \Map_{\calD^{\op}}(d^\vee,d'^\vee)\simeq \Map_{\calD}(d'^\vee,d^\vee)\simeq \Map_{\calD}(\mathbbm{1},d^\vee\otimes d')\,,
    \end{equation*}
    proving \ref{item:dualityadjunction}.
\end{proof}

\begin{lemma}
    \label{lem:laxtransformationapp}
    Let $\calC$ and $\calD$ be symmetric monoidal $\infty$-categories such that every object $d\in \calD$ is dualizable with dual $d^\vee$.
    Let $F:\calC\to \calD$ be symmetric monoidal and let $G:\calD\to \Sp$ be lax symmetric monoidal.
    Then there is a lax symmetric monoidal natural transformation of functors $\calC^{\op}\times \calC\to \Sp$
    \begin{equation*}
        \Sigma_+^\infty \Map_{\calC}(c,c')\longrightarrow G(F(c)^{\vee}\otimes F(c'))\,.
    \end{equation*}
    Its adjoint $\Map_{\calC}(c,c')\to \Omega^\infty G(F(c)^{\vee}\otimes F(c'))$ sends $f:c\to c'$ to the image under $G$ of the morphism $\mathbbm{1}\to F(c)^{\vee}\otimes F(c')$ corresponding to $F(f)$, precomposed with the unit $\mathbb{S}\to G(\mathbbm{1})$.
\end{lemma}
\begin{proof}
    By Lemma~\ref{lem:mappingspaces}, there are lax symmetric monoidal transformations of functors $\calC^{\op}\times\calC\to\Spc$
    \begin{equation*}
        \Map_{\calC}(c,c') \to \Map_{\calD}(F(c),F(c')) \simeq \Map_{\calD}(\mathbbm{1},F(c)^{\vee}\otimes F(c'))\,.
    \end{equation*}
    The functor $\Omega^\infty$ is lax symmetric monoidal and right adjoint to the symmetric monoidal functor $\Sigma_+^\infty$.
    The adjunction is monoidal.
    The functor $\Map_\calD(\mathbbm{1},-)$ is the unit of $\Fun(\calD,\Spc)$ with Day convolution and is therefore initial in $\CAlg(\Fun(\calD,\Spc))\simeq\Fun^{\mathrm{lax}}(\calD,\Spc)$~\cite[Cor.~3.2.1.9]{HA}.
    Thus, there is a unique lax symmetric monoidal transformation $\Map_\calD(\mathbbm{1},-)\to\Omega^\infty G$.
    Composing gives a lax symmetric monoidal transformation $\Map_\calC(c,c')\to\Omega^\infty G(F(c)^\vee\otimes F(c'))$.
    Its adjoint is the claimed transformation.
\end{proof}
\section{Comparison with existing orientations}
\label{app:comparison}
This section proves that the $\E_\infty$-orientation $M\Spin \to KO$ constructed in Theorem~\ref{thm:spin-orientation-proof}  is homotopic to
the $\E_\infty$-orientation of \cite{AndoHopkinsRezk}.
This is in turn an $\E_\infty$-refinement of the Atiyah--Bott--Shapiro orientation~\cite{ABS} and is homotopic to Joachim's model~\cite{JoachimCoherences}.
We will recall now the Hirzebruch characteristic series.

Let $R$ be a ring spectrum.
Let $V\to X$ be a vector bundle and let $u_\alpha,u_\beta\in R^n(\Th(V))$ be two Thom classes.
Then there is a \emph{difference class} $\delta(u_\alpha,u_\beta)\in R^0(X)^\times$ with $u_\alpha=\delta(u_\alpha,u_\beta)u_\beta$.
Let $s:X\to \Th(V)$ be the zero section.
The Euler classes $e_\alpha:=s^*u_\alpha$ and $e_\beta:=s^*u_\beta$ also satisfy $e_\alpha=\delta(u_\alpha,u_\beta)e_\beta$.

For any orientation $\beta:MSO\to R$, one can compare its rationalization to the canonical orientation
$\alpha:MSO\to H\Q \to R\otimes H\Q\simeq \pi_*R\otimes H\Q$ by comparing the Thom classes assigned to oriented vector bundles.
The universal complex line $L\to BU(1)$ carries an orientation.
Let $x=s^*u_\alpha\in H^2(BU(1);\Q)$ be the Euler class of $L$.
\begin{definition}[{\cite[Def.~3.10]{AndoHopkinsRezk}}]
    The \emph{Hirzebruch characteristic series} of $\beta$ is 
    \begin{equation*}
        K_\beta(x)=\delta(u_\alpha,u_\beta) = 1+o(x) \in H^0(BU(1);\pi_*R\otimes \Q)^\times \,.
    \end{equation*} 
\end{definition}

Note that $M\Spin \otimes H\Q\xrightarrow{\simeq} MSO\otimes H\Q$ and so rationalizing a $\Spin$-orientation gives rise to an $SO$-orientation with a well-defined characteristic series.

\begin{theorem}[{\cite[Thm.~6.1]{AndoHopkinsRezk}, \cite[Thm.~3.1.1]{NaumannSprang}}]
    \label{thm:ahr}
    The Hirzebruch characteristic series defines an injection
    \begin{equation*}
        \pi_0\Map_{\CAlg(\Sp)}(M\Spin, KO) \hookrightarrow 1+x^2\Q[[x^2]] \,.
    \end{equation*}    
\end{theorem}
\begin{remark}
    The precise statement in \cite{AndoHopkinsRezk} and \cite{NaumannSprang} uses the coefficients $b_{2k}$ of $K_\beta(x)=\exp\left(\sum_k 2b_{2k}(\beta)\frac{x^{2k}}{(2k)!}\right)$~\cite[Ex.~5.15]{AndoHopkinsRezk} and determines precise conditions to lie in the image.    
\end{remark}
Hence, two $\E_\infty$-maps $M\Spin \to KO$ are homotopic if and only if they have the same Hirzebruch series if and only if they give the same Thom or Euler classes.

We will now show that the $KO$-Euler class has a very simple description.
Let $X\to B\Spin$ be a spin vector bundle, equivalently encoded by a bimodule bundle $\Sigma:\Cl(V)\squigrightarrow \Cl_{\rk V}$.
The Thom class is an Euler class together with a trivialization over $S(V)$, which we can ignore for now.
Via \eqref{eq:deflaxmonoidaltrans}, Theorem~\ref{thm:spin-orientation-proof}, constructs a natural transformation, natural in $x\in X$:
\begin{equation*}
    \begin{tikzcd}
        *\ar[r,"{[\Cl(V_x)]}"]& k^{C_2}_{\Cl(V_x)} \ar[r,"\Sigma_x"] & k^{C_2}_{\Cl_{\rk V}} \to \Sigma^{\rk V}KO  \,,
    \end{tikzcd}
\end{equation*}
where the composite assigns to every fiber the graded $\Cl_{\rk V}$-module $\Sigma_x$.
That is, the Euler class is the bundle $[\Sigma]\in K^{\gr}_{\Cl_n}(X)\cong KO^n(X)$.
\begin{proposition}
    \label{prop:comparison}
    The Hirzebruch characteristic series of the $\E_\infty$-orientation $\beta\colon M\Spin\to KO$ of Theorem~\ref{thm:spinorientation} is
    \begin{equation*}
        K_\beta(x)=\frac{x}{e^{\frac{x}{2}}-e^{-\frac{x}{2}}}=\frac{x/2}{\sinh(x/2)} \,,
    \end{equation*}
    the characteristic series of the $\hat{A}$-genus.
    Consequently, $\beta$ is homotopic as an $\E_\infty$-ring map to the orientations of Ando--Hopkins--Rezk~\cite{AndoHopkinsRezk} and Joachim~\cite{JoachimCoherences}.
\end{proposition}
\begin{proof}
    Let $f\colon B\Spin(2)\to BSO(2)=BU(1)$ be induced by the double cover and let $W=f^*L$ be the universal spin bundle of rank~$2$.
    Under $\Spin(2)\cong U(1)$, let $H$ be the weight-one representation and $y=c_1(H)$, so that $H^{**}(B\Spin(2);\Q)=\Q[[y]]$, $W\cong H^{\otimes 2}$ and $f^*x=2y$.
    As $L$ is not spin, $K_\beta$ is defined through $M\Spin\otimes H\Q\simeq MSO\otimes H\Q$, and naturality of Thom classes gives $f^*K_\beta(x)=K_\beta(2y)$ as the difference class of the Thom classes of $W$.
    Restricting to the zero section yields $e_\alpha(W)=K_\beta(2y)\,e_\beta(W)$ with $e_\alpha(W)=2y$.
    Since $H^{**}(B\Spin(2);\pi_*KU\otimes\Q)$ has no zero divisors, it suffices to compute $e_\beta(W)$ rationally.

    The Euler class is $e_\beta(W)=[\Sigma_W]\in K^{\gr}_{\Cl_2}(B\Spin(2))\cong KO^2(B\Spin(2))$, where $\Sigma_W=E\Spin(2)\times_{\Spin(2)}\Cl_2$ and $\Spin(2)\subset\Cl_2^0$ acts by left multiplication.
    Difference classes are natural in the ring spectrum and $\pi_*KO\otimes\Q\to\pi_*KU\otimes\Q$ is injective, so we may complexify.
    The Morita equivalence $\Cxl_2\simeq\C$ given by the irreducible $\Cxl_2$-module $\Delta$ sends $[\Sigma_W\otimes\C]$ to $[\Cxl_2\otimes_{\Cxl_2}\Delta]=[\Delta]$, with $\Spin(2)$ acting through $\Cl_2^0$.
    As $e_1e_2$ acts by $\pm i$ on $\Delta^{0}$ and $\Delta^{1}$, we get $\Delta^0\cong H$ and $\Delta^1\cong H^{-1}$, hence $e_\beta(W)\mapsto[H]-[H^{-1}]\in KU^0(B\Spin(2))$.
    The rationalization $KU\to\prod_i\Sigma^{2i}H\Q$ is the Chern character, which maps this to $e^{y}-e^{-y}$.
    Thus $2y=K_\beta(2y)(e^{y}-e^{-y})$, and substituting $x=2y$ gives the claimed series.

    By \cite[Ex.~5.15]{AndoHopkinsRezk} this is the characteristic series of the Atiyah--Bott--Shapiro orientation, so Theorem~\ref{thm:ahr} shows that $\beta$ is homotopic to the orientation of \cite{AndoHopkinsRezk}, and hence to Joachim's~\cite{JoachimCoherences}.
\end{proof}

\bibliographystyle{alpha}
\bibliography{biblio.bib}{}

\end{document}

%% file: preamble.tex
\usepackage[utf8]{inputenc}
\usepackage[T1]{fontenc}

\usepackage{amsmath}
\usepackage{amsfonts}
\usepackage{amssymb}
\usepackage{amsthm}
\usepackage{mathtools}
\usepackage{bbm}

\usepackage{enumerate}
\usepackage{enumitem}
\setlist[enumerate]{before=\normalfont}
\usepackage{xcolor}
\usepackage{stmaryrd}
\usepackage{slashed}
\usepackage{setspace}
\usepackage{relsize}
\usepackage{comment}
\usepackage{wrapfig}

\usepackage{tikz-cd}
\usepackage{tikz}
\usetikzlibrary{decorations.markings}
\usetikzlibrary{decorations.pathmorphing}
\usetikzlibrary{cd}
\usepackage{tikz-3dplot}
\tikzset{snake it/.style={decorate, decoration=snake}}

\usepackage{tabularx}
\usepackage{booktabs}
\usepackage{array}
\usepackage{makecell}
\usepackage{geometry}
\usepackage{hyperref}

\newcommand{\Mod}{\operatorname{Mod}}

\newcommand{\Map}{\operatorname{Map}}
\newcommand{\map}{\operatorname{map}}
\newcommand{\sVect}{\operatorname{sVect}}
\newcommand{\Vect}{\operatorname{Vect}}

\newcommand{\rk}{\operatorname{rk}}

\newcommand{\Inn}{\mathcal{I}nn}

\newcommand{\Aut}{\operatorname{Aut}}
\newcommand{\Iso}{\operatorname{Iso}}
\newcommand{\sAut}{\operatorname{sAut}}
\newcommand{\sEnd}{\operatorname{sEnd}}
\newcommand{\Hom}{\operatorname{Hom}}
\newcommand{\Mor}{\operatorname{Mor}}

\newcommand{\End}{\operatorname{End}}

\DeclareMathOperator{\GL}{GL}

\newcommand{\Spin}{\operatorname{Spin}}
\newcommand{\Pin}{\operatorname{Pin}}
\newcommand{\id}{\operatorname{id}}

\newcommand{\ev}{\operatorname{ev}}
\newcommand{\coev}{\operatorname{coev}}

\newcommand{\sLine}{\operatorname{sLine}}

\newcommand{\Fun}{\operatorname{Fun}}
\newcommand{\TwArr}{\operatorname{TwArr}}
\newcommand{\Th}{\operatorname{Th}}
\DeclareMathOperator{\Lan}{Lan}

\newcommand{\Pic}{\operatorname{Pic}}

\newcommand{\op}{\operatorname{op}}
\newcommand{\slice}{\mathbin{\downarrow}}

\newcommand{\fib}{\operatorname{fib}}
\newcommand{\cofib}{\operatorname{cofib}}
\DeclareMathOperator{\Cl}{Cl}
\newcommand{\Cxl}{\mathbb{C}\mathrm{l}}

\newcommand{\Z}{\mathbb{Z}}
\newcommand{\N}{\mathbb{N}}
\newcommand{\Q}{\mathbb{Q}}
\newcommand{\C}{\mathbb{C}}
\newcommand{\R}{\mathbb{R}}
\newcommand{\E}{\mathbb{E}}

\newcommand{\calC}{\mathcal{C}}
\newcommand{\calD}{\mathcal{D}}
\newcommand{\calE}{\mathcal{E}}

\newcommand{\grotimes}{\hat{\otimes}}
\newcommand{\squigrightarrow}{\rightsquigarrow}
\newcommand{\xsquigrightarrow}[1]{\mathrel{\overset{#1}{\rightsquigarrow}}}

\newcommand{\Bim}{\operatorname{Bim}}
\newcommand{\const}{\mathrm{const}}
\newcommand{\gr}{\mathrm{gr}}

\providecommand{\sslash}{/\!\!/}

\newcommand{\Top}{\operatorname{Top}}

\newcommand{\sBan}{\mathrm{sBan}}

\newcommand{\Spc}{\operatorname{Spc}}

\newcommand{\Sp}{\operatorname{Sp}}
\DeclareMathOperator*{\colim}{colim}

\newcommand{\Grpd}{\operatorname{Grpd}}

\newcommand{\Alg}{\operatorname{Alg}}

\newcommand{\sSet}{\operatorname{sSet}}

\newcommand{\CAlg}{\mathsf{CAlg}}

\newcommand{\gp}{\mathrm{gp}}

\newtheorem{theorem}{Theorem}[section]
\newtheorem*{theorem*}{Theorem}
\newtheorem*{theoremA*}{Theorem~A}
\newtheorem*{theoremB*}{Theorem~B}
\newtheorem*{theoremC*}{Theorem~C}
\newtheorem{introtheorem}{Theorem}

\newtheorem{proposition}[theorem]{Proposition}
\newtheorem{lemma}[theorem]{Lemma}

\newtheorem{corollary}[theorem]{Corollary}

\theoremstyle{definition}

\newtheorem{definition}[theorem]{Definition}

\theoremstyle{remark}

\newtheorem{remark}[theorem]{Remark}

\newtheorem{example}[theorem]{Example}

\definecolor{Blue} {rgb} {0.282352,0.239215,0.803921}
\definecolor{Green} {rgb} {0.133333,0.545098,0.133333}

\newcounter{jfc}

\makeatletter
\newcommand{\nameditem}[1]{%
    \item[#1]
    \protected@edef\@currentlabel{#1}
}
\makeatother

\makeatletter
\renewcommand\paragraph{\@startsection{paragraph}{4}%
  \z@{.5\linespacing\@plus.7\linespacing}{-\fontdimen2\font}%
  {\normalfont\bfseries}}
\makeatother